\documentclass[10pt,reqno]{amsart}
\usepackage{amsfonts}
\usepackage{amsmath}
\usepackage{CJK}
\usepackage{amssymb}
\usepackage{latexsym,bm}
\usepackage[RGB]{xcolor}
\usepackage{mathrsfs}
\usepackage{amsthm}
\numberwithin{equation}{section}
\newtheorem{theorem}{Theorem}[section]

\newtheorem{lemma}[theorem]{Lemma}
\newtheorem{remark}{Remark}[section]

\newtheorem{proposition}[theorem]{Proposition}
\newcommand{\dif}{\mathrm{d}}
\newcommand{\di}{\mathrm{div}}

\newcommand{\SUM}[3]{\sum\limits_{{#1}={#2}}^{#3}}

\begin{document}
\title[Multi-d hyperbolized compressible Navier-Stokes equations]{Global existence versus blow-up for multi-d  hyperbolized compressible Navier-Stokes equations}
\author{Yuxi Hu and Reinhard Racke}
\thanks{\noindent Yuxi Hu, Department of Mathematics, China University of Mining and Technology, Beijing, 100083, P.R. China, yxhu86@163.com\\
\indent Reinhard Racke, Department of Mathematics and Statistics, University of
Konstanz, 78457 Konstanz, Germany, reinhard.racke@uni-konstanz.de}
\begin{abstract}
We consider the non-isentropic compressible Navier-Stokes equations  in two or three space dimensions for which the heat conduction of Fourier's law is replaced by Cattaneo's law and
the classical Newtonian flow is replaced by a revised Maxwell flow.  We show that a physical entropy exists for this new model. For two special cases, we show the global well-posedness of solutions with small initial data and the blow-up of solutions in finite time for a class of large initial data. Moreover, for vanishing relaxation parameters, the solutions (if it exists) are shown to converge to solutions of the classical system.
 {\bf AMS classification code}: 35 L 60, 35 B 44
\end{abstract}
\maketitle
\section{Introduction}
The classical compressible Navier-Stokes equations  in $\mathbb R^n\times \mathbb R^+$, $n=2,3$, are given by
\begin{align}\label{1.1}
\begin{cases}
\partial_t\rho+\mathrm{div} (\rho u)=0,\\
\partial_t(\rho u) +\mathrm{div} (\rho u\otimes u)+\nabla p=\mathrm{div}(S),\\
\rho\partial_t e+\rho u \cdot \nabla e + p\di u +\di q =S:\nabla u,\\
\end{cases}
\end{align}
with  the constitutive law for a Newtonian fluid,
\begin{align}\label{1.2}
S=\mu \left( \nabla u +\nabla u^T-\frac{2}{n} \mathrm{div}\, u I_n\right) +\lambda \mathrm{div} u I_n
\end{align}
and heat conduction given by Fourier's law,
\begin{align}\label{1.3}
q=-\kappa \nabla \theta.
\end{align}
Here, the functions $\rho, u, e, p, S, q, \theta$ denote the fluid density, velocity, specific internal energy per unit mass, pressure, stress tensor, heat flux and temperature, respectively.
$I_n$ denotes the identity matrix in $\mathbb R^n$. $\mu, \lambda$ are the shear and the
bulk viscosity constant, respectively. $\kappa$ is the constant heat conduction coefficient.

 By combining Newton'saw of viscosity and Hooke's law of elasticity, one obtains the following constitutive equation
\begin{align}\label{1.4}
\tau \dot{S}+S=\mu \left( \nabla u +\nabla u^T-\frac{2}{n} \mathrm{div} u I_n\right) +\lambda \, \mathrm{div}\, u I_n,
\end{align}
where $\dot {S}=\partial_t S+u\cdot \nabla S$, see \cite{YW14}. The positive parameter $\tau$ is the relaxation time describing the time lag in the response of the stress tensor to velocity gradient.
A  fluid  obeying equation \eqref{1.4} is called Maxwell flow, see also \cite{MA}. We should mention that,  even for a simple fluid, the effect of relaxation can not always be neglected, see \cite{PC} with the experiments of high-frequency (20GHZ) vibration of nanoscale mechanical devices immersed in water-glycerol mixtures. In 2014, Yong \cite{YW14} proposed a new model in which he divided  the stress tensor $S$ into $S_1+S_2 I_n$ and did relaxation for $S_1$ and $S_2$ in the following form
\begin{align}
\tau_2  \partial_t S_1+S_1=\mu (\nabla u+\nabla u^T-\frac{2}{n} \mathrm{div} u I_n),\label{1.5}\\
\tau_3  \partial_t S_2 +S_2=\lambda \mathrm{div} u. \label{1.6}
\end{align}
Note that there is no quadratic term $u\cdot \nabla S_1$ and $u\cdot \nabla S_2$  in \eqref{1.5} resp. \eqref{1.6}, and and thus does not have the property of Galilean invariance. The constitutive equations \eqref{1.5}-\eqref{1.6} are called revised Maxwell's law. A similar revised Maxwell model was considered by Chakraborty and Sader \cite{CS} for a compressible viscoelastic fluid, where they show that this model provides a general formalism with which to characterize the fluid-structure interaction of nanoscale mechanical devices vibrating in simple liquids.

On the other hand, although Fourier's law \eqref{1.3} plays an important role in experimental and applied physics, it has the drawbacks of an inherent infinite propagation speed of
signals. In order to overcome this paradox, Cattaneo \cite{CA} proposed the following constitutive equation for $q$,
\begin{align} \label{1.7}
\tau_1 \partial_t q+q+\kappa \nabla \theta=0,
\end{align}
which gives rise to heat waves with finite propagation speed. Here $\tau_1>0$ is the relaxation time.

In view of the above considerations, we investigate the following system in $\mathbb R^n \times \mathbb R_+$,
\begin{align} \label{1.8}
\begin{cases}
\partial_t\rho+\di (\rho u)=0,\\
\rho \partial_tu+\rho u \cdot\nabla u +\nabla p=\mu \di (\nabla u+\nabla u^T-\frac{2}{n} \di u I_n) +\nabla S_2,\\
\rho \partial_t e+\rho u\cdot \nabla e +p\di u +\di q =\mu (\nabla u+(\nabla u)^T-\frac{2}{n} \di u I_n): \nabla u+S_2\di u,\\
\tau_1 (\partial_t q+u \cdot \nabla q )+q+\kappa \nabla \theta=0,\\
\tau_3(\partial_tS_{2}+u\cdot \nabla S_2)+S_2=\lambda \di u,
\end{cases}
\end{align}
where we have taken $\tau_2=0$ in \eqref{1.5}.
That is, we do not consider a relaxation in $S_1$. This case is mathematically not yet accessible, even locally.
We consider  a Galilean invariant form of \eqref{1.6}. Besides, we consider the Galilean invariance form of Cattaneo's law $\eqref{1.8}_4$ which is proposed by Christov and Jordan \cite{CJo}.

We investigate the Cauchy problem to system \eqref{1.8} for the functions
\[
(\rho,u,\theta, q, S_2): \mathbb R^n\times[0,+\infty)\rightarrow \mathbb R_+\times \mathbb R^n\times \mathbb R_+ \times \mathbb R^n\times \mathbb R
\]
with initial condition
\begin{align}\label{1.9}
(\rho(x,0),u(x,0),\theta(x,0), q(x,0), S_2(x,0))=(\rho_0,u_0,\theta_0, q_0, S_{20}).
\end{align}
We are interested in the local and global well-posedness for small data, as well as in a blow-up for large data, both for the cases $\mu>0$ resp. $\mu=0$. The latter is not only motivated because it is mathematically accessible with respect to local existence and blow-up, but also with a physical background. In fact, there are recent studies determining the volume viscosity of a variety of gases which were found to be much larger (factor $10^4$) than the corresponding shear viscosities, see \cite{ShKu019}.

In order to close the system \eqref{1.8}, the equations for the thermodynamic variables $p$ and $e$ should be given. In this paper, we assume that
\begin{align}
e=C_v \theta+\frac{\tau_1}{\kappa \rho \theta} q^2 +\frac{\tau_3}{2\lambda \rho} S_2^2, \label{1.10}\\
p=R \rho \theta-\frac{\tau_1}{2\kappa \theta}q^2 -\frac{\tau_3}{2\lambda} S_2^2, \label{1.11}
\end{align}
with positive constants $C_v, R$ denoting the heat capacity at constant volume and the gas constant, respectively,
such that they satisfy the thermodynamic equation
$$\rho^2 e_\rho=p-\theta p_\theta.
$$

The dependence on $q^2$ term of the internal energy is indicated in paper \cite{CFO}, where they rigorously prove that such constitutive equations are consistent with
the second law of thermodynamics if and only if one use the relaxation equation $\eqref{1.8}_4$, see also \cite{CG, CHO, TA}. Since we also consider a relaxation for the stress tensor $S_2$, it is motivated, naturally,  by \cite{CFO} that the internal energy should also depend on $S_2$ in a quadratic form. Indeed, under the above constitutive laws, we show in Section 2 that there exists an entropy for the system \eqref{1.8}, which implies the compatibility with the second law of thermodynamics.

For classical compressible Navier-Stokes equations, the case $\tau_1=\tau_3=0$ in \eqref{1.8}, there are lots of results concerning the local and global well-posedness of strong  and/or weak solutions  because of its physical importance and mathematical challenges. In particular, the local existence and uniqueness of smooth solutions was established by Serrin \cite{SE} and Nash \cite{NA} for initial data far away from vacuum. Later, Matsumura and Nishida \cite{MN} got global smooth solutions for small initial data
without vacuum. For large data, Xin \cite{X}, Cho and Jin \cite{CJ} showed that smooth solutions must blow up in  finite time if the initial data has a vacuum state. See \cite{DA1,DA2, LI, JZ01, JZ03, FE} for global existence of weak solutions.

Neglecting the quadratic nonlinear terms $u\cdot \nabla q$ and $u \cdot \nabla S_2$ in $\eqref{1.8}_4$ resp. $\eqref{1.8}_5$, the cases $\tau_1>0, \tau_2=\tau_3=0$ (Cattaneo's law) and $\tau_1=0, \tau_2>0, \tau_3>0$ (revised Maxwell's law) have been studied in $\mathbb R^n,\, n\ge 2$, respectively, in our papers \cite{HR1, HR2}. For the \textit{one-dimensional} case, we had considered the relaxation both for $q$ and $S$ with Galilean invariance form. In \cite{HR3}, we showed the global existence of smooth solutions with small initial data and convergence to classical system as relaxation parameter goes to zero. In our paper with Wang \cite{HuRaWa022}, a blow-up result for large data was shown, hereby also yielding an interesting example, where the relaxed and the non-relaxed system are close to each other on the linearized level, globally for small data for the nonlinear system, and on any finite time horizon for the nonlinear one, but differ qualitatively for large data.

However, there are few results by considering both relaxations for $q$ and $S$ in the \textit{multi-dimensional} case. The aim of this paper is to solve this problem for two special cases with $\tau_1>0, \tau_2=0, \tau_3>0$ and $\mu>0$ resp. $\mu=0$.

Let $G:=\mathbb R^+\times \mathbb R^n\times \mathbb R \times \mathbb R^n\times \mathbb R$ denotes the physical state spaces of the unknowns $(\rho, u, \theta, q, S_2)$. Let $ G_0$ and $G_1$ be any convex compact subset of $G$ such that $G_0\subset\subset G_1 \subset\subset G$. The following theorem is considered under the assumption $\mu>0$. We have global existence for small data.

Defining the energy term
\begin{align}\label{3.4}
E(t):=\sup_{0\le \tau \le t} \|(\rho-1,u,\theta-1,q,S_2)(\tau,\cdot)\|_{H^3}^2+\int_0^t \left(\|(\nabla \rho, \nabla \theta)\|_{H^2}^2 + \|(q,S_2)\|_{H^3}^2 +\|\nabla u \|_{H^3}^2 \right)\dif t
\end{align}
and
$$
E_0:=E(0),
$$
we will prove first a global existence theorem for small data,
\begin{theorem}\label{th1.1}
Let $\tau_1>0, \tau_2=0, \tau_3>0$ and $\mu>0$. Suppose for the initial data $$(\rho_0-1, u_0, \theta_0-1, q_0, S_{20})\in H^3.$$ Then, there exists a small constant $\delta>0$ such that
if $E_0<\delta$, then the initial value problem \eqref{1.8}, \eqref{1.9} has a unique solution $(\rho, u, \theta, q, S_2)$ globally in time such that
$(\rho-1, u, \theta-1, q, S_2)\in C(0,+\infty; H^3)$, $(\nabla \rho, \nabla \theta)\in L^2(0, +\infty; H^2)$, $\nabla u \in L^2 (0, +\infty; H^3)$, $(q, S_2)\in L^2(0, +\infty; H^3)$.

For any $t>0$ we have
\begin{align}\label{1.12}
\|(\rho-1, u, \theta-1, q, S_2)\|_{H^3}^2 +\int_0^t \left( \|(\nabla \rho ,\nabla \theta)\|_{H^2}^2+\|\nabla u\|_{H^3}^2+\|(q, S_2)\|_{H^3}^2 \right)\dif t \le C E_0^2,
\end{align}
where $C$ is a constant being independent of $t$ and of the initial data. Moreover, the solution decays in the sense
\begin{align}\label{1.13}
\|\nabla (\rho, u, \theta, q ,S_2)\|_{L^2}\rightarrow 0 \quad \mathrm{as}\, t\rightarrow \infty.
\end{align}
\end{theorem}

For $\mu=0$,  the system \eqref{1.8} is reduced to a purely hyperbolic one with zero-order damping terms
\begin{align} \label{1.18}
\begin{cases}
\partial_t\rho+\di (\rho u)=0,\\
\rho \partial_tu+\rho u \cdot\nabla u +\nabla p=\nabla S_2,\\
\rho \partial_t e+\rho u\cdot \nabla e +p\di u +\di q =S_2\di u,\\
\tau_1 (\partial_t q+u \cdot \nabla q )+q+\kappa \nabla \theta=0,\\
\tau_3(\partial_tS_{2}+u\cdot \nabla S_2)+S_2=\lambda \di u.
\end{cases}
\end{align}

The existence of solutions to the system \eqref{1.18} with initial data \eqref{1.9}, even locally, is not immediately clear, since it is not symmetric nor strictly hyperbolic (nor hyperbolic-parabolic). By carefully calculating the eigenvalues and eigenvectors of the corresponding matrix in the associated first-order system, we show that \eqref{1.18} is a \textit{constantly hyperbolic} system, and thus has a local solution. Furthermore, we show  in the following theorem that smooth solutions can not exist globally for a class of large initial data.

We define some useful averaged quantities, as in \cite{HuRaWa022}:
\begin{align}
F(t):=\int_{\mathbb R^n}  x\cdot \rho(x,t) u(x,t)  \dif x, \label{5.9}\\
G(t):=\int_{\mathbb R^n} ({\mathcal E}(x,t)-\bar {\mathcal E}) \dif x, \label{5.10}
\end{align}
where ${\mathcal E}(x,t):=\rho(e+\frac{1}{2} u^2)$ is the total energy and $\bar {\mathcal E}:=\bar \rho(\bar e+\frac{1}{2} \bar u ^2)=C_v$.
\begin{theorem} \label{th1.3}
Let $(\rho, u, \theta, q, S_2)$ be the local solution to \eqref{1.18} on $[0, T_0)$ with initial data \eqref{1.9} (given by Theorem \ref{th5.1} below). Assume that the initial data $(\rho_0-1, u_0, \theta_0-1, q_0, S_{20})$ are compactly supported in
a ball $B_0(M)$ with radius $M>0$.  Moreover, we assume that
\begin{align}
G(0)>0, \label{1.16}\\
1<\gamma:=1+\frac{R}{C_v} <\frac{5}{3}
\end{align}
Then, there exists $u_0$ satisfying
\begin{align} \label{1.17}
F(0)> \max\left\{ \frac{128 \sigma \max \rho_0}{3(5-3\gamma)} , \frac{8 \sqrt{\pi  \max\rho_0}}{\sqrt{3(5-3\gamma)}}  \right\} M^4,
\end{align}
such that the length $T_0$ of the maximal interval of existence of a smooth solution $(\rho, u, \theta, q, S_2)$ of \eqref{1.18} is finite,
provided the compact support of the initial data is sufficiently large.
\end{theorem}

Summarizing we present first results for some multi-dimensional hyperbolic compressible Navier-Stokes equations including local well-posedness, global well-posedness for small data and blow-up for large data, using a variety of techniques to overcome the technical difficulties.

The paper is organized as follows. In Section 2 we derive an entropy equation for system \eqref{1.8} and present some preliminary lemmas. The global existence result Theorem \ref{th1.1} is
proved in Section~3. In Section 4 the local existence and the blow-up result for large data,  Theorem \ref{th1.3},
are proved. In Section 5 we give a remark on the singular limit as $\tau \rightarrow 0$.

We introduce some notation. $W^{m,p}=W^{m,p}(\mathbb R^n),\, 0\le m \le \infty,\, 1\le p \le \infty,$ denotes the usual Sobolev space with norm $\|\cdot\|_{W^{m,p}}$.
$H^m$ and $L^p$ stand for $W^{m,2} (\mathbb R^n)$ resp. $W^{0, p} (\mathbb R^n)$.
For  $m\times d$-matrices $B=(b_{jk})$, $M=(m_{jk})$, we denote $A:B=\sum_{j=1}^m\sum_{k=1}^d b_{jk}\,m_{jk}$, and $M^2=M:M$. For $m\in \mathbb{N}_0$ we denote by $\nabla^m\, v$ derivatives of $v$ of order $m$.

\section{Entropy equation and some preliminary inequalities }
In this part, we first derive an entropy equation for system \eqref{1.8}.  Defining the entropy
\begin{align} \label{2.1}
\eta:=C_v \ln \theta-R \ln \rho+\frac{\tau_1}{2\kappa \theta^2 \rho} q^2.
\end{align}
Similar to \cite{HR3}, we have for a local solution
\begin{lemma} \label{lem2.1}
\begin{align}\label{entropy}
\partial_t(\rho \eta)+\di (\rho u \eta)+\di \left(\frac{q}{\theta}\right)=\frac{q^2}{\kappa \theta^2}+\frac{S_2^2}{\theta \lambda} + \frac{ S_1^2}{2\mu\theta}.
\end{align}
\end{lemma}
\begin{proof}
Dividing equation $\eqref{1.8}_3$ by $\theta$, we have
\begin{align}
\frac{\rho}{\theta}\partial_t (C_v \theta+\frac{\tau_1}{\kappa \theta \rho} q^2 +\frac{\tau_3}{2\lambda \rho} S_2^2 )
+\frac{\rho }{\theta} u \cdot \nabla (C_v \theta +\frac{\tau_1}{\kappa \theta \rho} q^2+\frac{\tau_3}{2\lambda \rho} S_2^2)+R \rho \di u \qquad\qquad \nonumber\\
 -\frac{\tau_1}{2\kappa \theta^2} q^2 \di u  -\frac{\tau_3}{2\lambda \theta} S_2^2 \di u +\frac{\di q}{\theta}
=\frac{S_2 \di u}{\theta}+\frac{\mu}{\theta} (\nabla u+\nabla u^T-\frac{2}{n} \di u I_n): \nabla u.
\label{2.2a}
\end{align}
For the term $\frac{\rho}{\theta} \partial_t(\frac{\tau_1}{\kappa \theta {\rho}} q^2)$, we have
\begin{align}
\frac{\rho}{\theta}\partial_t (\frac{\tau_1}{\kappa \theta {\rho}} q^2)&=\rho \partial_t(\frac{\tau_1}{\kappa  {\rho}\theta^2} q^2)+\frac{\tau_1\partial_t\theta} {\kappa\theta^3}q^2 \nonumber\\
&=\rho \partial_t(\frac{\tau_1}{\kappa  {\rho}\theta^2} q^2)-\frac{1}{2}  \partial_t(\frac{1}{\theta^2}) \frac{\tau_1}{\kappa} q^2 \nonumber\\
&=\rho \partial_t(\frac{\tau_1}{\kappa  {\rho}\theta^2} q^2)-\frac{1}{2}  \partial_t(\frac{\tau_1}{\kappa \theta^2 } q^2)+  \frac{\tau_1 q \partial_tq}{\kappa \theta^2}\nonumber\\
&=\rho \partial_t(\frac{\tau_1}{\kappa  {\rho}\theta^2} q^2)-\frac{1}{2}  \partial_t(\frac{\tau_1}{\kappa \theta^2 } q^2) -  \frac{\tau_1 q}{\kappa \theta^2} u\cdot \nabla q- \frac{q^2}{\kappa \theta^2}- \frac{q \cdot \nabla \theta}{\theta^2}. \label{2.2}
\end{align}
For the term $\frac{\rho }{\theta} u \cdot \nabla (\frac{\tau_1}{\kappa \theta {\rho}} q^2)$, we get
\begin{align}
\frac{\rho }{\theta}  u \cdot \nabla (\frac{\tau_1}{\kappa \theta {\rho}} q^2)&=\rho u \cdot \nabla (\frac{\tau_1 }{\kappa  {\rho}\theta^2} q^2) + u \frac{\nabla \theta}{\theta^3} \frac{\tau_1}{\kappa} q^2 \nonumber\\
&=\rho u \cdot \nabla (\frac{\tau_1 }{\kappa  {\rho} \theta^2} q^2)-  u \cdot \nabla (\frac{1}{2} \frac{\tau_1}{\kappa \theta^2} q^2)
+  \frac{\tau_1q}{\kappa \theta^2} u \cdot \nabla q.  \label{2.3}
\end{align}
Then we have
\begin{align}
&\frac{\rho}{\theta} \partial_t(\frac{\tau_1}{\kappa \theta {\rho}} q^2)+\frac{\rho }{\theta} u \cdot \nabla (\frac{\tau_1}{\kappa \theta  {\rho}} q^2)
-\frac{\tau_1}{2\kappa \theta^2} q^2 \di u\nonumber\\
&=\left[\rho \partial_t(\frac{\tau_1}{\kappa  {\rho}\theta^2} q^2)+\rho u \cdot \nabla (\frac{\tau_1 }{\kappa  {\rho} \theta^2} q^2)\right]
-\left[\frac{1}{2} \partial_t(\frac{\tau_1}{\kappa \theta^2 } q^2)+  \di (u\frac{1}{2} \frac{\tau_1}{\kappa \theta^2} q^2)\right]
-\frac{q^2}{\kappa \theta^2}-  \frac{q \nabla \theta}{\theta^2} \nonumber\\
&=\left[ \partial_t(\frac{1}{2}\frac{\tau_1}{\kappa \theta^2 } q^2)+  \di (u\frac{1}{2} \frac{\tau_1}{\kappa \theta^2} q^2)\right]
-\frac{q^2}{\kappa \theta^2}-  \frac{ q \cdot \nabla \theta}{\theta^2}. \label{2.4}
\end{align}
On the other hand, we have
\begin{align}
&\frac{\rho}{\theta} \partial_t\left( \frac{\tau_3}{2\lambda\rho} S_2^2\right) +\frac{\rho }{\theta} u \cdot \nabla \left(\frac{\tau_3}{2\lambda \rho} S_2^2\right)
-\frac{\tau_3}{2\lambda \theta} S_2^2 \di u \nonumber \\
&=\frac{\tau_3}{\theta \lambda} S_2 \partial_tS_{2}-\frac{\tau_3 S_2^2}{2 \lambda \rho \theta} \partial_t\rho+\frac{\tau_3}{\theta \lambda} S_2 u \cdot \nabla S_2
-\frac{\tau_3 S_2^2}{2\lambda \rho \theta} u \nabla \rho
-\frac{\tau_3}{2\lambda \theta}S_2^2 \di u\nonumber\\
&=\frac{S_2}{\theta\lambda} \left( \tau_3(\partial_tS_{2}+u \cdot \nabla S_2)\right)-\frac{\tau_3 S_2^2}{2\lambda \rho \theta} (\partial_t\rho+u  \cdot \nabla \rho)
-\frac{\tau_3}{2\lambda \theta} S_2^2 \di u\nonumber\\
&=\frac{S_2}{\theta \lambda} (-S_2+\lambda \di u )+\frac{\tau_3 S_2^2}{2\lambda \rho \theta} \rho \di  u-
\frac{\tau_3}{2\lambda \theta}S_2^2  \di u\nonumber\\
&=-\frac{S_2^2}{\theta \lambda}+\frac{S_2\di u}{\theta}.  \label{2.5}
\end{align}
Moreover, notice that
\begin{align}\label{2.10}
(\nabla u+\nabla u^T-\frac{2}{n} \di u I_n ): \nabla u =\frac{1}{2} |\nabla u+\nabla u^T-\frac{2}{n} \di u I_n|^2 = \frac{S_1^2}{2\mu^2},
\end{align}
 since the matrix $ \nabla u +\nabla u^T-\frac{2}{n} \di u I_n$ is symmetric and traceless.
Therefore, we derive from \eqref{2.2a}, \eqref{2.4}, \eqref{2.5}, \eqref{2.10} the following equality
\begin{align}\label{2.6}
 \partial_t(C_v \rho  \ln \theta+\frac{\tau_1}{2\kappa \theta^2} q^2)+ \di (C_v\rho u  \ln \theta+u\frac{\tau_1}{2\kappa \theta^2} q^2)+R\rho \di u
 +\di (\frac{q}{\theta}) \nonumber\\
 =\frac{q^2}{\kappa \theta^2}+\frac{S_2^2}{\theta \lambda}+\frac{\mu}{2\theta} |\nabla u+\nabla u^T-\frac{2}{n} \di I_n|^2.
\end{align}
Now, we rewrite the equation $\eqref{1.8}_1$ as
\begin{align}\label{2.7}
R\rho \partial_t(\ln \rho )+R\rho u \nabla(\ln \rho)+R\rho \di u=0.
\end{align}
Combining the equations \eqref{2.6} and \eqref{2.7}, we get the desired result.
\end{proof}

\begin{remark} \label{re2.1}
In case $\tau_1=\tau_2=\tau_3=0$, the entropy $\eta$ defined in \eqref{2.1} is reduced to the quantity $C_v \ln \theta-R\ln \rho$ which is the entropy for classical
fluid dynamics. Moreover, the entropy equation \eqref{entropy} is reduced to the following classical entropy equation
\begin{align} \label{2.8}
\partial_t(\rho \eta) +\di (\rho u \eta) - \di (\frac{\kappa \nabla \theta}{\theta})=\frac{\kappa |\nabla \theta|^2}{\theta^2}+\frac{1}{\theta}\left( \frac{\mu}{2} |\nabla u+\nabla u^T-\frac{2}{n} \di u I_n|^2
+\lambda |\di u|^2 \right).
\end{align}
\end{remark}
Combining the equation of conservation of energy
\begin{align}\label{2.9}
\partial_t(\rho e+\frac{1}{2} \rho u^2)+\di (\rho u e +\frac{\rho u}{2} |u|^2)+\di (pu)+\di q=\di (S:u),
\end{align}
the mass equation $\eqref{1.8}_1$ and the entropy equation \eqref{entropy}, we get
\begin{lemma} \label{lem-dissenteq}
\begin{align}
\partial_t \left[C_v \rho(\theta-\ln \theta -1)+R (\rho\ln \rho-\rho+1) +(1-\frac{1}{2\theta}) \frac{\tau_1}{\kappa \theta} q^2
 +\frac{1}{2}\rho u^2 +\frac{\tau_3}{2 \lambda} S_2^2\right] \nonumber \\
+ \di \left[C_v \rho u (\theta-\ln \theta -1)+u(1-\frac{1}{2\theta})\frac{\tau_1}{\kappa \theta} q^2+\frac{\tau_3}{2\lambda} u S_2^2
+R \rho u \ln \rho-R\rho u-\frac{q}{\theta}+\frac{1}{2}\rho u |u|^2\right. \nonumber \\
\left.+pu+q-\mu u(\nabla u+\nabla u^T-\frac{2}{n} \di u I_n)-S_2 u\right]+\frac{q^2}{\kappa \theta^2}+\frac{S_2^2}{\theta \lambda}+\frac{\mu}{2\theta}|\nabla u+\nabla u^T-\frac{2}{n} \di u I_n|^2=0.
\label{dissenteq}
\end{align}
\end{lemma}
Denoting $V:=(\rho,u,\theta,q,S_2)$, and writing \eqref{dissenteq} as
$$
\partial_t\eta_1(V) + \di [\zeta(V)] +\frac{q^2}{\kappa \theta^2}+\frac{S_2^2}{\theta \lambda} +\frac{S_1^2}{2\mu\theta}=0,
$$
 we have a dissipative relation for the convex entropy pair $(\eta_1,\zeta)$ with convex $\eta_1$.

This equation will later on imply in particular the lower energy estimates of the solutions which are crucial to get the global existence of solutions for system \eqref{1.8}.
Next, we present some inequalities which is frequently used in the proof of our main results.
\begin{lemma} \label{le2.2}
In space dimensions $n=2,3$ we have the standard Sobolev imbeddings

i) $H^2 \hookrightarrow L^\infty$.\\

ii) $H^1 \hookrightarrow L^p$, \quad for $2\le p \le 6$.
\end{lemma}
The following Moser-type inequalities will be used in subsequent  sections and can be found as a standard tool for example in \cite{Ma84, Ra015}.
\begin{lemma} \label{moser}
(i) Let $r,m,n \in \mathbb N, \, 1 < p \leq \infty, \, h
\in C^r(\mathbb R^m), \, B \, := \, \| h \|_{C^r(\overline{B(0,1)})}.$
Then there is a constant $C=C(r,m,n,p) > 0$ such that for all
$w=(w_1,\dots, w_m) \in W^{r,p} (\mathbb R^n) \cap L^\infty (\mathbb R^n)$ with $\| w \|_{L^\infty}
\leq 1$ the inequality
\begin{equation}\label{mos1}
\| \nabla^r h(w) \|_{L^p} \leq C \, B \, \| \nabla^r w \|_{L^p}
\end{equation}
holds.\\[1em]
(ii) Let $m \in \mathbb N.$ Then there is a constant
$C=C(m,n) > 0$ such that for all $f,g \in W^{m,2} \cap L^\infty$ and $\alpha
\in \mathbb N^n_0, | \alpha | \leq m$, the following inequalities hold:
\begin{eqnarray}
\| \nabla^{|\alpha|} (fg) \|_2 &\leq &C( \| f \|_{L^\infty} \| \nabla^mg \|_2
 +  \| \nabla^mf \|_2 \| g \|_{L^\infty} ),\\ \mbox{}\nonumber\\
\| \nabla^{|\alpha|} (fg) - f \nabla^{|\alpha|} g \|_2 &\leq &C( \| \nabla f
\|_{L^\infty} \| \nabla^{m-1} g \|_2  +  \| \nabla^mf \|_2 \| g \|_{L^\infty} ). \\ \nonumber
\end{eqnarray}
\end{lemma}

\section{Global existence for small data: proof of Theorem \ref{th1.1}}
We first derive the following equation for $\theta$:
\begin{align}\label{3.1}
\rho e_\theta \partial_t\theta+ (\rho u e_\theta-\frac{2q}{\theta}) \nabla \theta +\theta p_\theta \di u +\di q=\frac{2}{\kappa \theta} q^2+\frac{1}{\lambda} S_2^2
+\frac{\mu}{2} |\nabla u+\nabla u^T-\frac{2}{n} \di u I_n|^2 .
\end{align}
\textit{Proof of (\eqref{3.1}.} Since
\begin{align*}
\rho \partial_te = \rho e_\theta \partial_t\theta +\rho e_\rho \partial_t\rho +\rho e_q \partial_tq+\rho e_{S_2} \partial_tS_{2},\\
\rho u \cdot \nabla e= \rho e_\theta u \cdot \nabla \theta +\rho e_\rho u \cdot \nabla \rho +\rho e_q  u \cdot \nabla q  +\rho e_{S_2} u\cdot \nabla S_2,
\end{align*}
we have
\begin{align*}
&\rho \partial_te+\rho u \nabla e \\
&= \rho e_\theta \partial_t\theta+\rho e_\theta u\cdot \nabla \theta +\rho e_\rho (\partial_t\rho+u \cdot \nabla \rho) + \rho e_q(\partial_tq+u \cdot \nabla q)+\rho e_{S_2} (\partial_tS_{2}+u\cdot \nabla S_2)\\
&=\rho e_\theta \partial_t\theta+\rho e_\theta u\cdot \nabla \theta-\rho^2 e_\rho \di u +\rho e_q(-\frac{1}{\tau_1} q-\frac{\kappa}{\tau_1} \nabla \theta)+\rho e_{S_2} (-\frac{1}{\tau_3} S_2+\frac{1}{\tau_3} \lambda \di u)\\
&=\rho e_\theta \partial_t\theta+\rho e_\theta u\cdot \nabla \theta-\rho^2 e_\rho \di u +\rho \frac{2\tau_1}{\kappa \rho \theta} q(-\frac{1}{\tau_1} q-\frac{\kappa}{\tau_1} \nabla \theta)
+\rho \frac{\tau_3}{\lambda \rho} S_2(-\frac{1}{\tau_3} S_2+\frac{1}{\tau_3} \lambda \di u)\\
&=\rho e_\theta \partial_t\theta+(\rho e_\theta u-\frac{2q}{\theta} )\cdot \nabla \theta-\rho^2 e_\rho \di u -\frac{2}{\kappa \theta}q^2 -\frac{1}{\lambda} S_2^2+S_2 \di u.
\end{align*}
By the valid thermodynamic equation $ \rho^2 e_\rho=p-\theta p_\theta$ and \eqref{2.10}, we derive the desired equation \eqref{3.1}. \hfill $\Box$\\
So, we have
\begin{align}\label{3.2}
\begin{cases}
\partial_t\rho+\di (\rho u)=0,\\
\rho \partial_tu+\rho u \cdot\nabla u +\nabla p=\mu(\triangle u+\frac{n-2}{n} \nabla \di u ) +\nabla S_2,\\
\rho e_\theta \partial_t\theta+ (\rho u e_\theta-\frac{2q}{\theta}) \nabla \theta +\theta p_\theta \di u +\di q=\frac{2}{\kappa \theta} q^2
+\frac{1}{\lambda} S_2^2+\frac{\mu}{2} |\nabla u+(\nabla u)^T-\frac{2}{n} \di u I_n|^2, \\
\tau_1 (\partial_t q+u \cdot \nabla q )+q+\kappa \nabla \theta=0,\\
\tau_3(\partial_tS_{2}+u\cdot \nabla S_2)+S_2=\lambda \di u.
\end{cases}
\end{align}
We first show the \textit{local} existence of solutions in the following theorem.
\begin{theorem} \label{th3.1}
Suppose that the initial data $(\rho_0-1, u_0, \theta_0-1, q_0, S_{20})\in H^3$ and $(\rho_0, u_0, \theta_0, q_0, S_{20})\in G_0$. Then, for each convex open subset $G_1$ satisfying
$G_0\subset \subset G_1 \subset \subset G$, there exists $T=T\left( \|(\rho_0,\dots,S_{20})\|_{H^3} \right)>0$, such that the system \eqref{3.2} has an unique classical solution $(\rho, u, \theta, q, S_2)$ satisfying
\begin{align*}
(\rho-1, \theta-1, q, S_2) \in C([0,T], H^3)\cap C^1([0,T], H^2)\\
u\in C([0,T], H^3)\cap C^1([0,T], H^1)
\end{align*}
and
\begin{align*}
(\rho, u, \theta, q, S_2)\in G_1,\quad \forall (x,t)\in \mathbb R^3 \times [0,T].
\end{align*}
\end{theorem}
\begin{proof}
First, we write the above system as a hyperbolic-parabolic system:
\begin{align*}
\begin{cases}
A_0^1 \partial_t V + \SUM j 1 n A_j \partial_{x_j} V =f(V,\nabla u),\\
A_0^2 \partial_t u -\SUM i 1 n \SUM j 1 n B_{ij} \partial_{x_ix_j} u=g(V,u,\nabla V, \nabla u),
\end{cases}
\end{align*}
where $V=(\rho, \theta, q, S_2)$ and
\begin{align*}
A_0^1=\left(
 \begin{array}{cccc}
 1&0&0&0\\
 0&\rho e_\theta& 0&0\\
 0&0&\frac{\tau_1}{\kappa}&0\\
 0&0&0& \tau_3
 \end{array}
 \right),\quad
 \SUM j 1 nA_j\xi_j= \left(
 \begin{array}{cccc}
 u\cdot \xi & 0&0&0\\
 0& \rho e_\theta u\cdot \xi & \xi^T&0\\
 0&\xi&u\cdot \xi &0\\
 0&0&0& u \cdot \xi
 \end{array}
 \right),\\
 f(V,\nabla u)= \left(
 \begin{array}{c}
 \rho \di u\\
 -\theta p_\theta \di u +\frac{2}{\kappa \theta} q^2
+\frac{1}{\lambda} S_2^2+\frac{\mu}{2} |\nabla u+(\nabla u)^T-\frac{2}{n} \di u I_n|^2\\
-\frac{q}{\kappa}\\
\lambda \di u - S_2
\end{array}
\right),\\
 A_0^2= \rho I_n,\quad \SUM i 1 n \SUM j 1 nB_{ij}\xi_i\xi_j =\mu (I_n+\frac{n-2}{n}\xi \xi^T),\\
  g(V,u,\nabla V,\nabla u)=-\rho u\cdot \nabla u-p_\rho \nabla \rho -p_\theta \nabla \theta -p_q \nabla q-(p_{S_2}-1) \nabla S_2,
\end{align*}
for $\xi \in \mathbb{R}^n$ with $|\xi|=1$. We observe that $A_0^1, A_0^2$ are symmetric and positive definite matrices, $A_j$ is a symmetric matrix, $B_{ij}$ are symmetric matrices, and $\SUM i 1 n \SUM j 1 n B_{ij}\xi_i\xi_j$ is positive definite for all $\xi \in \mathbb S^{n-1}$. For the equilibrium state $\bar V =(1,1,0,0), \bar u=0$, we have $f(\bar V,0)=0$ and
$g(\bar V, \bar u, 0,0)=0$. Then the local existence theorem follows immediately by using the results of Kawashima, see \cite{KW} or the appendix in \cite{Ra015}.
\end{proof}

For the global existence for small data, we do \emph{not} use the results of Kawashima \cite{KW} by checking the so-called Kawashima condition. Instead, we use energy estimates and the entropy equation to get the desired result. Since we have quadratic terms in the equation $\eqref{3.2}_3$, the Sobolev regularity in Kawashima \cite{KW} in $\mathbb{R}^3$ is at least $H^4$. Here we obtain $H^3$-estimates; Kawashima's results do not apply to our problem.
\begin{proposition}\label{pro3.1}
There is a constant $\delta>0$ such that if $E(t)\le \delta,\, t\in [0,T]$, then for $t\in [0,T]$
\begin{align}\label{3.5}
E(t)\le C \|(\rho_0-1,u_0,\theta_0-1, q_0, S_0\|_{H^3}^2+C E(t)^\frac{3}{2},
\end{align}
where $C$ is a constant being independent of $\delta$ and of the initial data.
\end{proposition}

The global existence of solutions for small data follows from Proposition \ref{pro3.1} and Theorem \ref{th3.1} immediately by continuing a local solution, if the initial data are chosen to be sufficiently small.

The following series of Lemmas is devoted to prove  Proposition \ref{pro3.1}.
First, from Lemma \ref{lem-dissenteq}, we derive the following lower energy estimates for a given local solution, cp. \cite{HR3} for similar arguments.
\begin{lemma} \label{le3.1}
\begin{align*}
\frac{\dif }{\dif t} \int_{\mathbb R^3}  \left(C_v \rho(\theta-\ln \theta -1)+R (\rho\ln \rho-\rho+1) +(1-\frac{1}{2\theta}) \frac{\tau_1}{\kappa \theta} q^2
+\frac{1}{2}\rho u^2 +\frac{\tau_3}{2 \lambda} S_2^2\right) \dif x \\
+ \int_{\mathbb R^3} \left(\frac{q^2}{\kappa \theta^2}+\frac{S_2^2}{\theta \lambda}+\frac{\mu}{2\theta} |\nabla u+(\nabla u)^T-\frac{2}{n} \di uI_n|^2\right) \dif x =0.
\end{align*}
Moreover, there exists  $\delta_0>0$ such that if $E(t)\le \delta_0$ for $t\in[0,T]$,  we have
\begin{align}\label{3.6}
\|(\rho-1, u, \theta-1, q, S_2)(\cdot,t)\|_{L^\infty} \le \frac{1}{4}
\end{align}
and
\begin{align}\label{3.7}
\int_{\mathbb R^3} \left( (\rho-1)^2+(\theta-1)^2 +u^2+q^2+S_2^2 \right)(x,t) \dif x +\int_0^t \int_{\mathbb R^3} \left( q^2+S_2^2+|\nabla u|^2+|\di u|^2 \right) \dif x \dif t \le C E_0.
\end{align}
\end{lemma}
 In the sequel we will assume
$$
E(t)\le \delta_0, \qquad \mbox{ for } t\in[0,T].
$$
Before we concentrate on the higher-order a priori estimates, by using the equations \eqref{3.2}, we can prove
\begin{lemma} \label{le3.2}
There is a small constant $\delta>0$ such that if $E(t)\le \delta,\, t\in [0,T]$, then for $t\in [0,T]$
\begin{align}\label{3.8}
\|(\partial_t\rho, \partial_t\theta, \partial_tq, \partial_tS_2)\|_{H^2}^2+\|\partial_tu\|_{H^1}^2+ \int_0^t \|(\partial_t \rho, \partial_t \theta, \partial_t q, \partial_t S_2, \partial_t u)\|_{H^2}^2 \dif t \le C E(t).
\end{align}
\end{lemma}
\begin{proof}
We only give a short proof for the estimates of $\partial_tu$. The other terms can be estimated in a similar way.  From equation $\eqref{3.2}_2$, we have
\begin{align*}
\|\partial_tu\|_{H^1}&\le \|u\cdot \nabla u+\frac{p_\rho}{\rho} \nabla \rho+\frac{p_\theta}{\rho}  \nabla \theta +\frac{p_q}{\rho} \nabla q+\frac{p_{S_2}}{\rho}  \nabla S_2+\frac{1}{\rho} \nabla S_2+\frac{\mu}{\rho} (\triangle u+\frac{n-2}{n} \nabla \di u)\|_{H^1}\\
&\le C E(t)^\frac{1}{2}+\|\frac{\mu}{\rho} (\triangle u+\frac{n-2}{n} \nabla \di u)\|_{H^1}.
\end{align*}
Since
\begin{align*}
\|\frac{\mu}{\rho} (\triangle u+\frac{n-2}{n} \nabla \di u)\|_{L^2}
\le C \|u\|_{H^2}\le C E(t)^\frac{1}{2}
\end{align*}
and
\begin{align*}
&\| \nabla \left( \frac{\mu}{\rho} (\triangle u+\frac{n-2}{n} \nabla \di u )\right)\|_{L^2}\\
&\le C \left( \|\frac{1}{\rho^2} \nabla \rho\|_{L^3} \|\triangle u+\frac{n-2}{n} \nabla \di u \|_{L^6}+ \|\frac{1}{\rho}\|_{L^\infty} \|\triangle u+\frac{n-2}{n} \nabla \di u\|_{H^1}\right)\\
&\le C( \|\nabla \rho\|_{H^1} \|u\|_{H^3}+\|u\|_{H^3})\le C E(t)^\frac{1}{2}.
\end{align*}
Thus, we have
\begin{align*}
\|\partial_t u\|_{H^1}^2 \le C E(t).
\end{align*}
On the other hand, using similar methods as above, we have
\begin{align*}
\|\partial_t u\|_{H^2}\le C \|(\nabla \rho, \nabla u, \nabla \theta, \nabla q, \nabla S_2, \triangle u)\|_{H^2}
\end{align*}
which implies by the definition of $E(t)$ that
\begin{align*}
\int_0^t \|\partial_t u\|_{H^2}^2 \dif t \le C E(t).
\end{align*}

\end{proof}

Now, for $1\le |\alpha| \le 3$, taking $\nabla\nabla^{|\alpha|-1}$ to the equations $\eqref{3.2}_1, \eqref{3.2}_3, \eqref{3.2}_4, \eqref{3.2}_5$ and $\nabla^ {|\alpha|-1}$ to the equation $\eqref{3.2}_2$,  we derive the following system
\begin{align}\label{high}
\begin{cases}
&\partial_t \nabla \nabla^{|\alpha|-1} \rho+ u\cdot \nabla (\nabla \nabla^{|\alpha|-1} \rho) +\rho \nabla \di \nabla^{|\alpha|-1} u= -[ \nabla \nabla^{|\alpha|-1} (u\cdot \nabla \rho)-u\cdot \nabla(\nabla \nabla ^{|\alpha|-1} \rho)]\\
&\quad\qquad \qquad\qquad \qquad \qquad\qquad\qquad \qquad \qquad  -[\nabla \nabla^{|\alpha|-1}(\rho \di u)-\rho \nabla \di \nabla^{|\alpha|-1} u]=:f_1,\\
&\rho \partial_t \nabla^{|\alpha|-1} u-\mu( \triangle \nabla^{|\alpha|-1} u+\frac{n-2}{n} \nabla \di \nabla^{|\alpha|-1} u) +p_\rho\nabla \nabla^{|\alpha|-1} \rho+p_\theta \nabla \nabla^{|\alpha|-1} \theta -\nabla\nabla^{|\alpha|-1} S_2\\
&+\nabla^{|\alpha|-1} (p_q \nabla q+p_{S_2} \nabla S_2)
=-[\nabla^{|\alpha|-1} (p_\rho \nabla \rho)-p_\rho \nabla \nabla^{|\alpha|-1} \rho]-[\nabla^{|\alpha|-1}(p_\theta\nabla\theta)-p_\theta \nabla \nabla^{|\alpha|-1} \theta]\\
&\qquad \qquad\qquad \qquad\qquad \qquad \qquad -[\nabla^{|\alpha|-1} (\rho \partial_t u)-\rho \partial_t \nabla^{|\alpha|-1} u]
 -\nabla^{|\alpha|-1} (\rho u\cdot \nabla u)=:f_2,\\
&\rho e_\theta \partial_t \nabla \nabla^{|\alpha|-1} \theta+(\rho u e_\theta-\frac{2q}{\theta})\cdot \nabla(\nabla \nabla^{|\alpha|-1} \theta)+\theta p_\theta \nabla \di \nabla^{|\alpha|-1} u+ \nabla\di \nabla^{|\alpha|-1} q\\
&=-(\nabla \nabla^{|\alpha|-1}(\rho e_\theta \partial_t\theta)-\rho e_\theta \partial_t  \nabla \nabla^{|\alpha|-1} \theta)-[\nabla \nabla^{|\alpha|-1}(\theta p_\theta \di u)-\theta p_\theta  \nabla \di \nabla^{|\alpha|-1} u]\\
&+ \nabla \nabla^{|\alpha|-1}\left( \frac{2}{\kappa \theta} q^2+ \frac{1}{\lambda} S_2^2 \right)
-\left(\nabla \nabla^{|\alpha|-1} ((\rho u e_\theta-\frac{2q}{\theta})\nabla \theta)-(\rho u e_\theta-\frac{2q}{\theta})\cdot \nabla(\nabla \nabla^{|\alpha|-1} \theta)\right)\\
&\qquad \qquad\qquad \qquad\qquad \qquad\qquad \qquad \qquad+\frac{\mu}{2} \nabla \nabla^{|\alpha|-1} \left( |(\nabla u+\nabla u)^T-\frac{2}{n} \di u I_n|^2\right)=:f_3,\\
&\tau_1 (\partial_t \nabla \nabla^{|\alpha|-1} q +u \cdot \nabla (\nabla \nabla^{|\alpha|-1} q) +\nabla \nabla^{|\alpha|-1} q +\kappa  \nabla^2\nabla^{|\alpha|-1} \theta\\
&\qquad \qquad\qquad \qquad\qquad \qquad\qquad \qquad \qquad=-\tau_1 (\nabla\nabla^{|\alpha|-1}(u\cdot \nabla q)-u\cdot \nabla(\nabla \nabla^{|\alpha|-1} q)=:f_4,\\
&\tau_3(\partial_t \nabla \nabla^{|\alpha|-1} S_2+u \cdot \nabla(\nabla \nabla^{|\alpha|-1} S_2)+ \nabla\nabla^{|\alpha|-1} S_2-\lambda \nabla  \di\nabla^{|\alpha|-1}  u\\
&\qquad\qquad \qquad\qquad \qquad\qquad \qquad \qquad=-\tau_3(\nabla\nabla^{|\alpha|-1}(u\cdot \nabla S_2)-u\cdot \nabla (\nabla \nabla^{|\alpha|-1} S_2)=:f_5.
\end{cases}
\end{align}
The following lemma gives the estimates of the right-hand side of the above system.
\begin{lemma} \label{le3.3}
For any $1\le |\alpha| \le 3$, we have
\begin{align*}
\|f_1\|_{L^2} \le C (E(t))^\frac{1}{2} \|(\nabla \rho, \nabla u)\|_{H^2},\\
\|f_2\|_{L^2}\le C (E(t))^\frac{1}{2} ( \|(\nabla \rho, \nabla \theta)\|_{H^1}+ \| \nabla u\|_{H^2}+\|\partial_t u\|_{H^1}),\\
\|f_3\|_{L^2} \le C (E(t))^\frac{1}{2}( \|(\nabla \rho, \nabla u, \nabla \theta, \nabla q, \nabla S_2)\|_{H^2}+\|\nabla u\|_{H^3}+\|\partial_t \theta\|_{H^2}),\\
\|f_4\|_{L^2} \le C (E(t))^\frac{1}{2} \|(\nabla q, \nabla u)\|_{H^2},\qquad \|f_5\|_{L^2} \le C (E(t))^\frac{1}{2}  \|(\nabla u, \nabla S_2)\|_{H^2}.
\end{align*}
\end{lemma}
\begin{proof}
By Sobolev's imbedding theorem and Moser-type inequalities, we have
\begin{align*}
\|f_1\|_{L^2} &\le \| \nabla \nabla^{|\alpha|-1}(u\cdot \nabla \rho)- u \cdot \nabla (\nabla \nabla^{|\alpha|-1}) \rho\|_{L^2}+\|\nabla\nabla^{|\alpha|-1} (\rho \di u)-\rho \nabla \di \nabla^{|\alpha|-1} u\|_{L^2}\\
&\le C\left( \|\nabla u\|_{L^\infty} \|\nabla^{|\alpha|} \rho\|_{L^2}+ \|\nabla \rho\|_{L^\infty} \|\nabla^{|\alpha|} u\|_{L^2}\right)\\
&\le C (E(t))^\frac{1}{2} \|(\nabla^{|\alpha|} \rho, \nabla^{|\alpha|} u)\|_{L^2}\le C (E(t))^\frac{1}{2} \|(\nabla \rho, \nabla u)\|_{H^2} .
\end{align*}
Noting that $p_\rho=R \theta, p_\theta=R\rho+\frac{\tau_1}{2\kappa \theta^2} q^2$, we get for $|\alpha|=1$ that
\begin{align*}
\nabla^{|\alpha|-1} (p_\rho \nabla \rho)-p_\rho \nabla \nabla^{|\alpha|-1}\rho=0
\end{align*}
and for $|\alpha|\ge 2$  that
\begin{align*}
\| \nabla^{|\alpha|-1} (p_\rho \nabla \rho)-p_\rho \nabla \nabla^{|\alpha|-1}\rho\|_{L^2}& \le C \|(\nabla \rho, \nabla \theta)\|_{L^\infty} \|(\nabla^{|\alpha|-1} \rho, \nabla^{|\alpha|-1} \theta)\|_{L^2}\\
&\le C (E(t))^\frac{1}{2} \|(\nabla \rho, \nabla \theta)\|_{H^2}.
\end{align*}
Similarly, we have for $|\alpha| \ge 1$ that
\begin{align*}
\|\nabla^{|\alpha|-1} (p_\theta \nabla \theta)- p_\theta \nabla \nabla^{|\alpha|-1} \theta \| \le C (E(t))^\frac{1}{2} \|(\nabla \rho, \nabla \theta, \nabla q)\|_{H^2}.
\end{align*}
Now we estimate the term $\| \nabla^{|\alpha|-1}(\rho \partial_t u)-\rho \partial_t \nabla^{|\alpha|-1} u \|_{L^2}$. For $|\alpha|=1$, this term vanishes.
For $|\alpha|=2$, we have as a typical term
\begin{align*}
\|\partial_{x_i}(\rho \partial_t u)- \rho \partial_t \partial_{x_i}u \|_{L^2}=\|(\partial_{x_i} \rho) \partial_t u\|_{L^2} \le (E(t))^\frac{1}{2} \|\partial_t u\|_{L^2}.
\end{align*}
For $|\alpha|=3$, we have as a typical term
\begin{align*}
&\| \partial_{x_i} \partial_{x_j} (\rho \partial_t u) -\rho \partial_t \partial_{x_i} \partial_{x_j} u\|_{L^2}\\
=& \|(\partial_{x_i} \partial_{x_j} \rho) \partial_t u+ (\partial_{x_j} \rho) \partial_t \partial_{x_i} u +\partial_{x_i} \rho \partial_{x_j} \partial_t u\|_{L^2}\\
\le& C( \| \nabla^2 \rho\|_{L^6} \|\partial_t u \|_{L^3} +\|\nabla \rho\|_{L^\infty} \|\partial_t u\|_{H^1}\\
 \le& C \|\nabla \rho\|_{H^2} \|\partial_t u\|_{H^1} \le C (E(t))^\frac{1}{2} \|\partial_t u\|_{H^1}.
\end{align*}
So, for $1 \le |\alpha| \le 3$, we get
\begin{align*}
\| \nabla^{|\alpha|-1}(\rho \partial_t u)-\rho \partial_t \nabla^{|\alpha|-1} u \|_{L^2} \le C (E(t))^\frac{1}{2} \|\partial_t u\|_{H^1}.
\end{align*}
On the other hand,
\begin{align*}
\|\nabla^{|\alpha|-1}(\rho u \cdot \nabla u)\|_{L^2} &\le C \|\rho u\|_{L^\infty} \|\nabla^{|\alpha|} u\|_{L^2}+ C\|\nabla u\|_{L^\infty} \| \nabla^{|\alpha|-1}(\rho u)\|_{L^2}\\
&\le C (E(t))^\frac{1}{2} ( \|\nabla u\|_{H^2}+ \|\nabla \rho\|_{H^1}).
\end{align*}
Combining the above estimates, we get the estimate for $f_2$.  The terms $f_3, f_4, f_5$ can be estimated in a similar way, we omit the details. This finishes the proof of Lemma \ref{le3.3}.
\end{proof}
The next lemma gives the higher-order a priori estimates of the solutions.
\begin{lemma}\label{le3.4}
For any $0\le t\le T$, we have
\begin{align} \label{3.9}
\|(\nabla \rho, \nabla u, \nabla \theta, \nabla q, \nabla S_2)\|_{H^2}^2 +\int_0^t (\|(\nabla q, \nabla S_2)\|_{H^2}^2+\| \nabla u\|_{H^3}^2 ) \dif t \le C (E_0+E(t)^\frac{3}{2}).
\end{align}
\end{lemma}
\begin{proof}
Multiply $\eqref{high}_1, \eqref{high}_2, \eqref{high}_3, \eqref{high}_4, \eqref{high}_5$ for $|\alpha|\geq 2$ by $\frac{p_\rho}{\rho}\nabla\nabla^{|\alpha|-1} \rho$, $-\triangle \nabla^{|\alpha|-1} u$, $\frac{1}{\theta}\nabla\nabla^{|\alpha|-1} \theta$, $\frac{1}{\kappa \theta} \nabla \nabla^{|\alpha|-1} q$, $\frac{1}{\lambda}\nabla \nabla^{|\alpha|-1} S_2$, respectively, -- for $|\alpha|=1$ take the multipliers  $\frac{p_\rho}{\rho}\nabla  \rho$, $-\triangle  u$, $\frac{1}{\theta}\nabla  \theta$, $\frac{1}{\kappa \theta} \nabla   q$, $\frac{1}{\lambda}\nabla  S_2$, respectively --, and summing up the results, we get
\begin{align}
&\frac{\dif}{\dif t} \int \left[ \frac{p_\rho}{2\rho} |\nabla \nabla^{|\alpha|-1} \rho|^2+\frac{\rho}{2}  |\nabla \nabla^{|\alpha|-1} u|^2 +\frac{ \rho e_\theta}{2\theta} |\nabla \nabla ^{|\alpha|-1} \theta|^2 +\frac{\tau_1}{2\kappa\theta} |\nabla \nabla^{|\alpha|-1} q|^2+\frac{\tau_3}{2\lambda} |\nabla \nabla^{|\alpha|-1} S_2|^2 \right] \dif x \nonumber \\
&+\int \left(\frac{1}{\kappa \theta} |\nabla \nabla^{|\alpha|-1} q|^2+\frac{1}{\lambda}|\nabla \nabla^{|\alpha|-1} S_2|^2  \right)\dif x+\mu \int ( \triangle \nabla^{|\alpha|-1} u+\frac{n-2}{n} \nabla \di \nabla^{|\alpha|-1} u ) \triangle \nabla^{|\alpha|-1} u \dif x \nonumber \\
&=\SUM i 1 4 G_i+\SUM i 1 4 H_i+\SUM i 1 5 F_i+D+L \label{New1}
\end{align}
where
$$
\int \dots \dif x := \int_{\mathbb{R}^n} \dots \dif x,
$$
and
\begin{align*}
G_1:=\int u \cdot \nabla(\nabla \nabla^{|\alpha|-1} \rho) \cdot \frac{p_\rho}{\rho} \nabla \nabla^{|\alpha|-1} \rho \dif x, \quad
G_2:=\int (\rho u e_\theta-\frac{2q}{\theta})\cdot \nabla(\nabla \nabla^{|\alpha|-1} \theta) \cdot \frac{1}{\theta} \nabla\nabla^{|\alpha|-1} \theta \dif x,\\
G_3:=\int \tau_1 u \cdot \nabla(\nabla \nabla^{|\alpha|-1} q) \cdot \frac{1}{\kappa \theta} \nabla \nabla^{|\alpha|-1} q \dif x,\quad
G_4:= \int \tau_3 u \cdot \nabla (\nabla \nabla^{|\alpha|-1} S_2)\cdot \frac{1}{\lambda} \nabla \nabla^{|\alpha|-1} S_2 \dif x,
\end{align*}
\begin{align*}
H_1:= \int \left((\rho \nabla \di \nabla^{|\alpha|-1} u) \cdot \frac{p_\rho}{\rho} \nabla \nabla^{|\alpha|-1} \rho - p_\rho \nabla \nabla^{|\alpha|-1} \rho \cdot \triangle \nabla^{|\alpha|-1} u\right)\dif x,\\
H_2:=\int \left(- p_\theta \nabla \nabla^{|\alpha|-1} \theta \cdot \triangle \nabla^{|\alpha|-1} u+\theta p_\theta \nabla \di \nabla^{|\alpha|-1} u\cdot \frac{1}{\theta} \nabla \nabla^{|\alpha|-1} \theta \right) \dif x,\\
H_3:=\int \left( \frac{1}{\theta} \nabla \di \nabla^{|\alpha|-1} q \cdot \nabla \nabla^{|\alpha|-1} \theta +\kappa \nabla^2 \nabla^{|\alpha|-1} \theta \cdot \frac{1}{\kappa \theta} \nabla\nabla^{|\alpha|-1} q\right)\dif x,\\
H_4=\int \left(\nabla \nabla^{|\alpha|-1} S_2\cdot \triangle  \nabla^{|\alpha|-1}u-\lambda\nabla \di \nabla^{|\alpha|-1} u \cdot \frac{1}{\lambda} \nabla \nabla^{|\alpha|-1} S_2\right)\dif x,
\end{align*}
\begin{align*}
F_1:=\int f_1 \cdot \nabla \nabla^{|\alpha|-1} \rho \dif x,\quad F_2:=-\int f_2\cdot \triangle \nabla^{|\alpha|-1} u\dif x, \quad F_3:=\int f_3 \cdot \nabla \nabla^{|\alpha|-1} \theta \dif x,\\
F_4:=\int f_4\cdot \frac{1}{\kappa} \nabla \nabla^{|\alpha|-1} q \dif x,\quad F_5:= \int f_5 \cdot \nabla \nabla^{|\alpha|-1} S_2 \dif x,
\end{align*}
\begin{align*}
&D:=\int \frac{1}{2} \left(\partial_t \left( \frac{p_\rho}{\rho} \right)|\nabla \nabla^{|\alpha|-1} \rho|^2 +\frac{1}{2} \partial_t \rho |\nabla \nabla^{|\alpha|-1} u|^2  \right.\\
&\qquad \qquad \qquad\qquad \qquad\qquad \qquad \qquad \left.+\frac{1}{2} \partial_t \left( \frac{\rho e_\theta}{\theta}\right) |\nabla \nabla^{|\alpha|-1} \theta|^2 - \nabla \rho \cdot \partial_t \nabla^{|\alpha|-1} u \cdot \nabla \nabla^{|\alpha|-1} u \right)\dif x,
\end{align*}
\begin{align*}
L:=\int -\nabla^{|\alpha|-1} (p_q \nabla q +p_{S_2} \nabla S_2 ) \cdot \triangle \nabla^{|\alpha|-1} u \dif x.
\end{align*}
Now we estimate the right hand side of the equation \eqref{New1}.
\begin{align*}
G_1&=\int-\di (\frac{p_\rho}{\rho} u) \cdot \frac{1}{2} |\nabla \nabla^{|\alpha|-1} \rho|^2 \dif x \\
&\le C \|\nabla (\rho, \theta, u)\|_{L^\infty} \int |\nabla \nabla^{|\alpha|-1} \rho|^2 \dif x \le C (E(t))^\frac{1}{2} \|\nabla \rho \|_{H^2}^2\dif x,
\end{align*}
while $G_i,\, 2\le i\le 4$ can be estimated similarly.  So, we get
\begin{align}\label{3.10}
\SUM i 1 4 G_i \le C (E(t))^\frac{1}{2} \|(\nabla \rho, \nabla \theta, \nabla q, \nabla S_2)\|_{H^2}^2.
\end{align}
 Integrating by part, we get
\begin{align*}
H_1&=\int p_\rho \nabla \nabla^{|\alpha|-1} \rho \cdot ( \nabla \di \nabla^{|\alpha|-1} u- \triangle \nabla ^{|\alpha|-1} u) \dif x \\
&= \int p_\rho \nabla \nabla^{|\alpha|-1} \rho \cdot \nabla \times \nabla \times (\nabla^{|\alpha|-1} u) \dif x \\
&=\int (\nabla p_\rho \times \nabla \nabla^{|\alpha|-1} \rho) \cdot \nabla \times (\nabla^{|\alpha|-1} u) \dif x \le C E(t)^\frac{1}{2} \| \nabla \rho\|_{H^2} \|\nabla u\|_{H^2}.
\end{align*}
Here we used the relation
$$
\Delta u=  \nabla \di u-\nabla \times \nabla \times u,
$$
which also holds in two space dimensions. There,   for a scalar field $f$ the rotation is given as the vector $\nabla\times f:=(\partial_{x_2} f, - \partial_{x_1} f)^T$, resp. for a vector field $F$ as the scalar $\nabla \times F := \partial_{x_1}F_2 - \partial_{x_2}F_1$.
Similarly, we have
\begin{align*}
H_2&=\int p_\theta \nabla \nabla^{|\alpha|-1} \theta ( \nabla \di \nabla^{|\alpha|-1} u- \triangle \nabla^{|\alpha|-1} u) \dif x \\
&=\int p_\theta \nabla \nabla^{|\alpha|-1} \theta \cdot \nabla \times \nabla \times \nabla^{|\alpha|-1} u \dif x \le C E(t)^\frac{1}{2} \|\nabla \theta\|_{H^2} \|\nabla u\|_{H^2}.
\end{align*}
On the other hand, we have
\begin{align*}
H_3&=\int \frac{1}{\theta} \left( \nabla \di \nabla^{|\alpha|-1} q \cdot \nabla \nabla^{|\alpha|-1} \theta +\nabla^2\nabla^{|\alpha|-1} \theta \cdot \nabla \nabla^{|\alpha|-1} q \right) \dif x \\
&=\int \frac{1}{\theta} \nabla \nabla^{|\alpha|-1} \theta \cdot ( \nabla \di \nabla^{|\alpha|-1} q- \triangle \nabla^{|\alpha|-1} q) \\
&=\int \nabla  \left(\frac{1}{\theta}\right)\cdot  \nabla \nabla^{|\alpha|-1} \theta \cdot \nabla \times \nabla^{|\alpha|-1} q \dif x  \le C E(t)^\frac{1}{2} \|\nabla \theta\|_{H^2} \|\nabla q\|_{H^2}.
\end{align*}

Moreover,
\begin{align*}
H_4&=\int \nabla^{|\alpha|-1} \nabla S_2(\triangle \nabla^{|\alpha|-1} u-\nabla \di \nabla^{|\alpha|-1} u)\dif x \\
&=-\int \nabla^{|\alpha|-1}\nabla S_2\cdot \nabla \times \nabla \times \nabla^{|\alpha|-1} u \dif x=0.
\end{align*}

 So, we derive that
\begin{align} \label{3.11}
\SUM i 1 4 H_i \le  C (E(t)^\frac{1}{2})  \|(\nabla \rho, \nabla u, \nabla \theta, \nabla q)\|_{H^2}^2.
\end{align}
Using Lemma \ref{le3.3}, we get
\begin{align} \label{3.12}
\SUM i 1 5 F_i \le C E(t)^\frac{1}{2} \left( \|(\nabla \rho, \nabla u, \nabla \theta, \nabla q, \nabla S_2)\|_{H^2}^2 +\|\nabla u\|_{H^3}^2 +\|\partial_t u\|_{H^1}^2+\|\partial_t \theta\|_{H^1}^2 \right).
\end{align}
For the last terms $D$ and $L$ on the right-hand side of \eqref{New1},  we have
\begin{align*}
D& \le \|(\rho_t, \theta_t, q_t, (S_2)_t, \nabla \rho)\|_{L^\infty} \int \left(|\nabla \nabla ^{|\alpha|-1} \rho|^2 +|\nabla \nabla^{|\alpha|-1} u|^2
+|\nabla \nabla^{|\alpha|-1} \theta|^2+|\nabla^{|\alpha|-1} u_t|^2\right)\dif x\\
&\le C E(t)^\frac{1}{2} \left( \|(\nabla \rho, \nabla u, \nabla \theta, \nabla q)\|_{H^2}^2+ \|u_t\|_{H^2}^2 \right) \dif x.
\end{align*}
To estimate $L$, we first take $|\alpha|=1$, typically
 \begin{align*}
L&=\int (p_q \nabla q+p_{S_2} \nabla S_2) \triangle u \dif x \\
&\le C \|(p_q, p_{S_2})\|_{L^\infty} \int( |\nabla q|^2+|\nabla S_{2}|^2+|\triangle u|^2 )\dif x \\
&\le C E(t)^\frac{1}{2}(\|(\nabla q, \nabla S_2)\|_{H^2}^2+\|u\|_{H^2}^2),
\end{align*}
where we have used the fact that $p_q=-\frac{\tau_1}{\kappa \theta} q, p_{S_2}=-\frac{\tau_3}{\lambda} S_2$,
while for $2\le |\alpha| \le 3$,
\begin{align*}
L&=-\int \nabla^{|\alpha|-1} (p_q \nabla q+p_{S_2} \nabla S_2) \triangle \nabla^{|\alpha|-1} u \dif x \\
&\le C \|(p_q, \nabla q, p_{S_2}, \nabla S_2)\|_{L^\infty}\int \left( |\nabla \nabla^{|\alpha|-1} q|^2 +|\nabla^{|\alpha|-1} q |^2+ |\nabla^{|\alpha|-1} \theta|^2+|\nabla \nabla^{|\alpha|-1} S_2|^2\right.\\
&\left.\qquad +|\nabla^{|\alpha|-1}S_2|^2+|\triangle \nabla^{|\alpha|-1} u|^2 \right)\dif x \\
&\le C E(t)^\frac{1}{2}(\|(\nabla q, \nabla \theta, \nabla S_2)\|_{H^2}^2+\|\nabla u\|_{H^3}^2).
\end{align*}
For the last term on the left hand side of \eqref{New1}, we observe that
\begin{align*}
&\int \nabla \di \nabla^{|\alpha|-1} u \cdot \triangle \nabla^{|\alpha|-1} u \dif x\\
=&\int \nabla \di \nabla^{|\alpha|-1} u \cdot (\nabla \di \nabla^{|\alpha|-1}u-\nabla\times \nabla \times \nabla^{|\alpha|-1} u)\dif x\\
=&\int \left| \nabla \di  \nabla^{|\alpha|-1}u\right|^2 \dif x,
\end{align*}
which gives
\begin{align}
\mu \int ( \triangle \nabla^{|\alpha|-1} u+\frac{n-2}{n} \nabla \di \nabla^{|\alpha|-1} u ) \triangle \nabla^{|\alpha|-1} u \dif x \nonumber\\
=\mu \int |\triangle \nabla^{|\alpha|-1} u|^2 +\frac{n-2}{n} \left| \nabla \di  \nabla^{|\alpha|-1}u\right|^2 \dif x \label{3.13}
\end{align}
Summing $|\alpha|$ from $1$ to $3$ and integrating the inequality \eqref{New1} over $(0,t)$,   using the elliptic inequality $\|u\|_{H^2}\le C \|\triangle u\|_{L^2}$ as well as Lemma \ref{le3.2}, we get the desired estimates \eqref{3.9}, which proves Lemma \ref{le3.4}.
 \end{proof}
 Now, we use the ideas from Kawashima \cite{KW} , see also \cite{UmKaSh84, SK85}, to estimate the missing last two terms: $\int_0^t \| \nabla \rho\|_{H^2}^2 \dif t$ and $\int_0^t \|\nabla \theta\|_{H^2}^2 \dif t$.

We first write the system \eqref{3.2} in symmetric hyperbolic-parabolic form.
Let $U:=(\rho, u, \theta, q, S_2)$ and $\bar U=(\bar \rho, \bar u, \bar \theta, \bar q, \bar S_2):=(1, 0, 1, 0, 0)$,  then we have
\begin{align} \label{3.14}
A^0(\bar U) U_t+\SUM j 1 n A^j(\bar U) U_{x_j}+\SUM j 1 n \SUM k 1 n B^{jk} (\bar U) U_{x_jx_k} +L(\bar U) U=F(U,\nabla U, \partial_t U),
\end{align}
where
\begin{align*}
A^0(\bar U)=\left(
\begin{array}{ccccc}
\bar p_\rho&0 &0 &0 &0\\
0& I_n&0 &0 &0\\
0 &0 & \bar e_\theta& 0&0\\
0 &0 &0 & \frac{\tau_1}{\kappa} &0\\
0 &0 &0 &0 & \frac{\tau_3}{\lambda}
 \end{array}\right),
 \SUM j 1 n  A^j(\bar U)\xi_j=\left(
 \begin{array}{ccccc}
 0& \bar p_\rho \xi^T&0&0&0\\
 \bar p_\rho \xi &0& \bar p_\theta \xi &0& -\xi\\
 0& \bar p_\theta \xi^T& 0 &\xi^T& 0\\
 0&0&\xi&0&0\\
 0& -\xi^T& 0&0&0
 \end{array}
 \right),\\
 \SUM j 1 n \SUM k 1 n B^{jk}\xi_j\xi_k=\left(
 \begin{array}{ccccc}
 0&0&0&0&0\\
 0& \mu I_n+ \frac{n-2}{n}\mu \xi \xi^T&0&0&0\\
 0&0&0&0&0\\
 0&0&0&0&0\\
 0&0&0&0&0\\
 \end{array}
 \right),
 L(\bar U)= \left(
 \begin{array}{ccccc}
  0&0&0&0&0\\
 0&0&0&0&0\\
 0&0&0&0&0\\
0&0&0&\frac{1}{\kappa}&0\\
 0&0&0&0&\frac{1}{\lambda}
 \end{array}
 \right),
 \end{align*}
$\xi=(\xi_1, \cdots, \xi_n)\in S^{n-1}$and  $F(U,\nabla U, U_t)\equiv (F_1, F_2, F_3, F_4, F_5)$ where
\begin{align*}
&F_1:=\bar p_\rho(-u\cdot  \nabla \rho-(\rho-1)\di u)\\
&F_2:= -(\rho-1) \partial_t u -(p_\rho-\bar p_\rho)\nabla \rho-(p_\theta-\bar p_\theta)\nabla \theta-p_q \nabla q-p_{S_2} \nabla S_2,\\
&F_3:=-(\rho e_\theta-\bar e_\theta)\partial_t \theta-(\rho u e_\theta-\frac{2q}{\theta})\nabla \theta-(\theta p_\theta-\bar p_\theta) \di u +\frac{2}{\kappa \theta}q^2+\frac{1}{\lambda} S_2^2
+\frac{\mu}{2} |\nabla u+ \nabla u^T-\frac{2}{n} \di u|^2,\\
&F_4:=-\frac{\tau_1}{\kappa} u\cdot \nabla q,\quad F_5:= -\frac{\tau_3}{\lambda} u\cdot \nabla S_2
\end{align*}
 $F$ is etimated in the following lemma.
\begin{lemma} \label{le3.5}
For $1 \le |\alpha| \le 3$, there exists a constant C such that
\begin{align*}
\|\nabla^{|\alpha|-1} F_1\|_{L^2}\le C E(t)^\frac{1}{2} \|(\nabla \rho, \nabla u)\|_{H^2},\\
\|\nabla^{|\alpha|-1} F_2\|_{L^2} \le C E(t)^\frac{1}{2} \left( \|(\nabla \rho, \nabla \theta, \nabla q, \nabla S_2)\|_{H^2}+\|\partial_t u\|_{H^2}\right),\\
\|\nabla^{|\alpha|-1} F_3\|_{L^2} \le C E(t)^\frac{1}{2} \left( \|(\nabla \theta, \partial_t \theta, \nabla u)\|_{H^2}+\|(q, S_2)\|_{H^2}\right),\\
\|\nabla ^{\alpha-1} F_4\|_{L^2} \le  C E(t)^\frac{1}{2} \left( \|\nabla q\|_{H^2}+\|\nabla u\|_{H^1}\right),\\
\|\nabla^{|\alpha|-1} F_5\|_{L^2} \le C E(t)^\frac{1}{2} \left( \|\nabla S_2\|_{H^2}+ \|\nabla u\|_{H^1}\right).
\end{align*}
\end{lemma}
\begin{proof}
For $|\alpha|=1$, we have
\begin{align*}
\|F_1\|_{L^2}\le \|(u, \rho-1)\|_{L^\infty} (\|\nabla \rho\|_{L^2}+\|\nabla u\|_{L^2})\le C E(t)^\frac{1}{2} (\|(\nabla \rho ,\nabla u)\|_{H^2}).
\end{align*}
For $2 \le |\alpha| \le 3$, we have
\begin{align*}
\|\nabla^{|\alpha|-1} F_1\|_{L^2} &\le C \|(u, \nabla \rho, \rho-1, \nabla u)\|_{L^\infty} (\|(\nabla^{|\alpha|} \rho, \nabla^{|\alpha|-1} u, \nabla^{|\alpha|} u, \nabla^{|\alpha|-1} \rho)\|_{L^2})\\
&\le C E(t)^\frac{1}{2} \|(\nabla \rho ,\nabla u)\|_{H^2}.
\end{align*}
The other terms in $F_i, \, 2\le i \le 5$ can be estimated in a similar way as above except the first term "$-(\rho-1)u_t$" in $F_2$ since $\|u_t\|_{L^\infty}$ is not bounded by $E(t)^\frac{1}{2}$.
In the following, we  only estimate this term.  Note that for $|\alpha|=1$,
\begin{align*}
\|(\rho-1)\partial_tu\|_{L^2} \le \|\rho-1\|_{L^\infty} \|\partial_tu\|_{L^2}\le C E(t)^\frac{1}{2} \|\partial_t u\|_{L^2}.
\end{align*}
For $|\alpha|=2$,
\begin{align*}
&\|\partial_{x_i} [(\rho-1)\partial_t u]\|_{L^2} \\
&=\| \partial_{x_i} \partial_t u+ (\rho-1) \partial_t\partial_{x_i} u\|\\
&\le C \|(\nabla \rho, \rho-1)\|_{L^\infty} \|\partial_t u\|_{H^1} \le C E(t)^\frac{1}{2} \|\partial_t u\|_{H^1}.
\end{align*}
For $|\alpha|=3$,
\begin{align*}
&\|\partial_{x_i}\partial_{x_j} [(\rho-1)\partial_t u ]\|_{L^2}\\
&=\| \partial_{x_i}\partial_{x_j} \rho \partial_t u+ \partial_{x_i} \partial_t \partial_{x_j} u+\partial_{x_j} \rho \partial_t \partial_{x_i} u+(\rho-1) \partial_t \partial_{x_i}\partial_{x_j} u\|\\
&\le \|\nabla^2\rho\|_{L^6} \|\partial_t u\|_{L^3}+\|\nabla \rho\|_{L^\infty} \|\nabla \partial_t u\|_{L^2}+\|\rho-1\|_{L^\infty} \|\partial_t u \|_{H^2}\\
&\le C E(t)^\frac{1}{2} \|\partial_t u\|_{H^2}.
\end{align*}
This finishes the proof of Lemma \ref{le3.5}.
\end{proof}
Now, we introduce a matrix $K^j$ for $j=1,2, \cdots, n$ as follows,
\begin{align} \label{3.15}
\SUM j 1 n K^j \xi_j=\epsilon \left(
\begin{array}{ccccc}
0& \bar p_\rho \xi^T&0&0&0\\
-\xi&0&0&0&0\\
0&0&0& \frac{\kappa N}{\tau_1} \xi^T&0\\
0&0& -\frac{N}{\bar e_\theta} \xi&0&\frac{\lambda}{\tau_3} \xi\\
0&0&0&-\frac{\kappa}{\tau_1} \xi^T&0
\end{array}
\right),
\end{align}
where $\epsilon>0$ and $N>0$ is yet arbitrary.
Then, we calculate
\begin{align} \label{3.16}
\SUM j 1 n K^j A^0(\bar U)\xi_j=\epsilon\left(
\begin{array}{ccccc}
0&\bar p_\rho \xi^T&0&0&0\\
-\bar p_\rho \xi&0&0&0&0\\
0&0&0& N\xi^T&0\\
0&0&-N\xi&0&\xi\\
0&0&0&-\xi^T&0
\end{array}\right)
\end{align}
which is an anti-symmetric matrix and
\begin{align*}
M:=\SUM j 1 n \SUM i 1 n \left(\frac{1}{2} \left( K^iA^j \xi_j\xi_i+(K^iA^j\xi_j\xi_i)^T\right)+ B^{ij}\xi_j\xi_i\right)+L\\
=\left(
\begin{array}{ccccc}
\epsilon (\bar p_\rho)^2 &0& \frac{\epsilon}{2} \bar p_\rho \bar p_\theta &0& -\frac{\epsilon}{2} \bar p_\rho\\
0& \mu I_n+(\frac{n-2}{n}\mu-\epsilon \bar p_\rho)\xi \xi^T&0&-\frac{\epsilon}{2}\left(\frac{\bar p_\theta N}{\bar e_\theta}+\frac{\lambda}{\tau_3}\right)\xi \xi^T&0\\
\frac{\epsilon}{2}\bar p_\rho \bar p_\theta&0&\epsilon \frac{\kappa N}{\tau_1}&0&- \epsilon \frac{\kappa}{2\tau_1}\\
0& -\frac{\epsilon}{2}\left(\frac{\bar p_\theta N}{\bar e_\theta}+\frac{\lambda}{\tau_3}\right)\xi \xi^T&0&\frac{1}{\kappa} I_n-\frac{\epsilon N}{\bar e_\theta} \xi \xi^T&0\\
-\frac{\epsilon}{2} \bar p_\rho&0&-\epsilon \frac{\kappa }{2\tau_1}&0&\frac{1}{\lambda}
\end{array}
\right)
\end{align*}
is symmetric and positive definite for some large $N$ and sufficiently small $\epsilon$. In fact, let $\eta=(\eta_1, \eta_2, \eta_3, \eta_4, \eta_5)^T\in \mathbb R \times \mathbb R^n\times \mathbb R\times \mathbb R^n \times \mathbb R$ be any vector, then
$$
\eta^T M \eta =(\eta_1, \eta_3, \eta_5)^T M_1 (\eta_1, \eta_3, \eta_5)+(\eta_2, \eta_4)^T M_2 (\eta_2, \eta_4),
$$
where
\begin{align*}
M_1=\left(
\begin{matrix}
\epsilon (\bar p_\rho)^2 & \frac{\epsilon}{2} \bar p_\rho \bar p_\theta & -\frac{\epsilon}{2} \bar p_\rho\\
\frac{\epsilon}{2}\bar p_\rho \bar p_\theta&\epsilon \frac{\kappa N}{\tau_1}&- \epsilon \frac{\kappa}{2\tau_1}\\
-\frac{\epsilon}{2} \bar p_\rho&-\epsilon \frac{\kappa }{2\tau_1}&\frac{1}{\lambda}
\end{matrix}
\right),
M_2=\left(
\begin{matrix}
 \mu I_n+(\frac{n-2}{n}\mu-\epsilon \bar p_\rho)\xi \xi^T&-\frac{\epsilon}{2}\left(\frac{\bar p_\theta N}{\bar e_\theta}+\frac{\lambda}{\tau_3}\right)\xi \xi^T\\
 -\frac{\epsilon}{2}\left(\frac{\bar p_\theta N}{\bar e_\theta}+\frac{\lambda}{\tau_3}\right)\xi \xi^T&\frac{1}{\kappa} I_n-\frac{\epsilon N}{\bar e_\theta} \xi \xi^T\\
\end{matrix}
\right)
\end{align*}
So, the positive definiteness of $M$ is equivalent to that of $M_1$ and $M_2$.  Let $d_j$ denote the $j$-th principle minor of the matrix $M_1$, then $d_1=\epsilon (\bar p_\rho)^2>0$ and
$$
d_2=\epsilon^2 (\bar p_\rho)^2\left(\frac{\kappa N}{\tau_1}-\frac{1}{4} (\bar p_\theta)^2 \right).
$$
We can choose $N$ independent of $\epsilon$ such that
$$
\frac{\kappa N}{\tau_1}> \frac{1}{4} (\bar p_\theta)^2
$$
implying $d_2>0$. Now, $N$ is fixed. For small $\epsilon$, we observe that
$$
d_3=\frac{1}{\lambda} (\bar p_\rho)^2\left(\frac{\kappa N}{\tau_1}-\frac{1}{4} (\bar p_\theta)^2\right) \epsilon^2+O(\epsilon^3), \quad \, \epsilon \rightarrow 0,
$$
which gives $d_3>0$ by choosing $\epsilon$ sufficiently small. Thus, $M_1$ is a positive definite matrix.  Furthermore, we observe that $M_2$ is a small symmetric perturbation of a positive definite matrix, which implies that $M_2$ is also positive definite for sufficiently small $\epsilon$. From now on $\epsilon$ is fixed.

Applying $K^i\nabla^{|\alpha|-1}$ to the system \eqref{3.14}, multiplying by $\nabla^{|\alpha|-1} U_{x_i}$ in $L^2$, and taking the sum over $i=1,\dots,n$, we obtain
\begin{align}
\SUM i 1 n < K^i A^0(\bar U) \nabla^{|\alpha|-1} U_t, \nabla^{|\alpha|-1} U_{x_i}>+\SUM i 1 n \SUM j 1 n < K^i A^j (\bar U)\nabla^{|\alpha|-1} U_{x_j}, \nabla^{|\alpha|-1} U_{x_i}>\nonumber\\
+<\SUM  i 1 n K^i \left( \SUM j 1 n \SUM k 1 n B^{jk}(\bar U) \nabla^{|\alpha|-1} U_{x_jx_k}+ L(\bar U)\nabla^{|\alpha|-1} U\right), \nabla^{|\alpha|-1} U_{x_i}>\nonumber\\
=\SUM i 1 n <K^i \nabla^{|\alpha|-1}F(U, \nabla U, U_t), \nabla^{|\alpha|-1} U_{x_i}>. \label{3.17}
\end{align}
Here $<\cdot, \cdot >$ means the $L^2$ inner product. Now, we estimate each term of \eqref{3.17} as follows.
\begin{align*}
&\SUM i 1 n < K^i A^0(\bar U) \nabla^{|\alpha|-1} U_t, \nabla^{|\alpha|-1} U_{x_i}>\\
=&\SUM i 1 n \frac{1}{2} \frac{\dif}{\dif t} < K^i A^0(\bar U) \nabla^{|\alpha|-1} U, \nabla^{|\alpha|-1} U_{x_i}>
\end{align*}
where the we use the fact $K^i A^0$ is an anti-symmetric matrix.
On the other hand, using the fact that $M$ defined above is a positive definite matrix, we derive that
\begin{align*}
&\SUM i 1 n \SUM j 1 n < K^i A^j (\bar U)\nabla^{|\alpha|-1} U_{x_j}, \nabla^{|\alpha|-1} U_{x_i}>\\
=&<\SUM i 1 n \SUM j 1 n\left(  \frac{1}{2} (K^i A^j(\bar U)+(K^i A^j(\bar U))^T)+B^{ij}(\bar U)\right)+L(\bar U) \nabla^{|\alpha|-1} U_{x_j}, \nabla^{|\alpha|-1} U_{x_i}>\\
&-<\SUM i 1 n \SUM j 1 n \left(B^{ij}(\bar U)+L(\bar U)\right) \nabla^{|\alpha|-1} U_{x_j}, \nabla^{|\alpha|-1} U_{x_i}>\\
&\ge \beta \|\nabla^{|\alpha|} U\|_{L^2}^2-C(\|(\nabla^{|\alpha|} u, \nabla^{|\alpha|} q, \nabla^{|\alpha|} S_2)\|_{L^2}^2,
\end{align*}
where $\beta=\beta(M)$ is some positive constant. The positivity of $M$ can be exploited using the Fourier transform.
Meanwhile, we have
\begin{align*}
&<\SUM  i 1 n K^i \left( \SUM j 1 n \SUM k 1 n B^{jk}(\bar U) \nabla^{|\alpha|-1} U_{x_jx_k}+ L(\bar U)\nabla^{|\alpha|-1} U\right), \nabla^{|\alpha|-1} U_{x_i}>\\
&\le C (\|\nabla u\|_{H^3}^2+ \|q\|_{H^3}^2+ \|S_2\|_{H^3}^2).
\end{align*}
Integrating the equation \eqref{3.17} over $(0,t)$ , combing the above estimates, using Lemmas \ref{le3.4} and \ref{le3.5}, we get the following result.
\begin{lemma} \label{le3.6}
There exists a constant $C>0$ such that
\begin{align} \label{3.18}
\int_0^t \|(\nabla \rho, \nabla \theta)\|_{H^2}^2 \dif t \le C(E_0+E(t)^\frac{3}{2}).
\end{align}
\end{lemma}

Combining Lemmas \ref{le3.4} and \ref{le3.6},  we get
\begin{lemma}\label{le3.7}
There exists a constant $C>0$ such that
\begin{align} \label{3.19}
\|(\nabla \rho,\nabla  u, \nabla  \theta, \nabla  q, \nabla  S_2)\|_{H^2}+ \int_0^t (\|(\nabla  \rho, \nabla  \theta, \nabla  q, \nabla  S_2)\|_{H^2}^2 +\|\nabla  u\|_{H^3}^2 ) \dif t \le C( E_0+ E(t)^{\frac{3}{2}}).
\end{align}
\end{lemma}
Thus, the proof of Proposition \ref{pro3.1} is finished by combining Lemmas \ref{le3.1} and \ref{le3.7}. Moreover, looking at the local existence Theorem \ref{th3.1} and the dependence on the norms of the initial data that determine the length of the existence interval, this allows to continue a local solution to a global one. Moreover, from inequality \eqref{1.12} and Lemma \ref{le3.2}, we have
\begin{align*}
\int_0^t \| \nabla (\rho, u, \theta, q, S_2)\|_{L^2} \dif t\le C
\end{align*}
and
\begin{align*}
\int_0^t \left| \frac{\dif}{\dif t} \|\nabla (\rho, u, \theta, q, S_2)\|_{L^2} \right|\dif t \le C,
\end{align*}
which implies the decay estimate \eqref{1.13} immediately and thus proves Theorem \ref{th1.1}. \hfill $\Box$


\section{Blow-up result: Proof of Theorem \ref{th1.3}}

 In this part, we consider the local existence and the blow-up of solutions for system \eqref{1.18}. To this end, we need the following assumption.\\
{\bf Assumption. } There exists $\delta>0$, sufficiently small, such that
\begin{align}
\min_{x\in \mathbb R^n}(\rho_0(x), \theta_0(x))>0,\quad \max_{x\in \mathbb R^n} (|q_0(x)|, |S_{20}(x)|)\le \frac{\delta}{2}. \label{5.3}
\end{align}
We remark that for the local solution then will hold
\begin{align}
\min_{(x,t) } (p_\rho, p_\theta, e_\theta)(t, x)>0\, \max_{(x,t) }|p_{S_2}(t,x)|<\frac{1}{2}. \label{5.4}
\end{align}
Replacing the equation $\eqref{1.18}_3$ by the equation for the temperature $\theta$ given in \eqref{3.1}, we have
\begin{align}\label{5.1}
\begin{cases}
\partial_t\rho+\di (\rho u)=0,\\
\rho \partial_tu+\rho u \cdot\nabla u +\nabla p=\nabla S_2,\\
\rho e_\theta \partial_t\theta+ (\rho u e_\theta-\frac{2q}{\theta}) \nabla \theta +\theta p_\theta \di u +\di q=\frac{2}{\kappa \theta} q^2
+\frac{1}{\lambda} S_2^2,\\
\tau_1 (\partial_t q+u \cdot \nabla q )+q+\kappa \nabla \theta=0,\\
\tau_3(\partial_tS_{2}+u\cdot \nabla S_2)+S_2=\lambda \di u.
\end{cases}
\end{align}
Now, we transform the above system into a first-order system for $V:=(\rho,u,\theta,q,S_2)^T$.
We have
\begin{align}\label{5.5}
\partial_tV+\SUM j 1 3 A_j(V) \partial_{x_j} V +B(V) V=F(V),
\end{align}
where
\begin{align*}
\SUM j 1 3 A_j \xi_j =\left(
\begin{array}{ccccc}
u\cdot \xi & \rho \xi^T& 0&0&0\\
\frac{p_\rho}{\rho} \xi & (u\cdot \xi) I_n& \frac{p_\theta}{\rho} \xi & \frac{p_q}{\rho} \xi^T& \frac{p_{S_2}-1}{\rho} \xi\\
0&\frac{\theta p_\theta}{\rho e_\theta} \xi^T & (u-\frac{2q}{\rho \theta e_\theta})\xi & \xi^T&0\\
0&0&\frac{\kappa}{\tau_1} \xi & (u\cdot \xi) I_n&0\\
0& -\frac{\lambda}{\tau_3} \xi^T & 0&0& u\cdot \xi
\end{array}
\right),\\
B(V)=diag\left\{0,0,0,\frac{1}{\tau_1},\frac{1}{\tau_3}\right\}, F(V)=\left(0,0, \frac{2}{\kappa \theta} q^2+\frac{1}{\lambda} S_2^2, 0,0\right)^T.
\end{align*}
Since the $(2n+3) \times (2n+3)$-matrix $\SUM j 1 n A_j \xi_j$ is \emph{not} symmetric, the system \eqref{5.5} is \emph{neither} symmetric-hyperbolic \emph{nor} strictly hyperbolic. So, the  local existence does not follow immediately by the classical theory of symmetric-hyperbolic or strictly hyperbolic systems. We shall show that the system \eqref{5.5} is a constantly hyperbolic system which also will imply a local existence theorem.

We look only at the three-dimensional case $n=3$, analogous arguments apply to the case $n=2$.
We first prove that the matrix $\SUM j 1 3 A_j\xi_j$ has nine linearly independent eigenvectors corresponding to five different eigenvalues (one eigenvalue is 5-fold).
For $|\xi|=1$ the characteristic polynomial for $\SUM j 1 3 A_j\xi_j$ is given by (also checked by \copyright Maple)
\begin{align} \label{5.6}
P(n,\Lambda,\xi):=\det\left(\SUM j 1 3 A_j\xi_j-\Lambda I_9\right) =(u\cdot \xi-\Lambda)^5 g(u\cdot \xi-\Lambda),
\end{align}
where
\begin{small}
\begin{align}
g(z):=&z^4-\frac{2q\cdot \xi}{\rho \theta e_\theta} z^3-\left( \frac{\kappa}{\tau_1}+\frac{\theta p_\theta^2}{\rho^2 e_\theta}+\frac{\lambda(1-p_{S_2})}{\rho \tau_3} +p_\rho\right) z^2 \nonumber \\
&+\left(\frac{\kappa \theta p_\theta}{\tau_1 \rho^2 e_\theta} p_q\cdot \xi+\left(\frac{\lambda (1-p_{S_2})}{\rho \tau_3}+p_\rho\right) \cdot \frac{2q\cdot \xi}{\rho \theta e_\theta} \right) z+\left(\frac{\lambda (1-p_{S_2})}{\rho \tau_3}+p_\rho\right) \cdot \frac{\kappa}{\tau_1}.  \label{5.7}
\end{align}
\end{small}
We observe the similarity of the characteristic polynomial $P(n,\Lambda,\xi)$ with the corresponding polynomial in one space dimension (n=1), given in \cite{HR3}, actually we have
\begin{equation}\label{charpol-n}
P(n,\Lambda,\xi)=(u\cdot \xi-\Lambda)^{2n-1} g(u\cdot\xi-\Lambda),\qquad \mbox{ for } n=1,2,3.
\end{equation}
A similar situation -- characteristic polynomial in space dimensions $n=2,3$ is given by the corresponding polynomial in space dimension $n=1$ times a power of linear polynomials -- is observed in linear thermoelasticity, see \cite{Ra015} or \cite{JR00}.

We first show that there exist five linearly independent eigenvectors corresponding to the eigenvalue $\Lambda=u\cdot\xi$. Let $W=(x_1,x_2,x_3,x_4,x_5,x_6,x_7,x_8,x_9)^T$ be an   eigenvector corresponding to the eigenvalue $\Lambda=u\cdot\xi$. Then we have
\begin{align*}
x_2 \xi_1+x_3 \xi_2+x_4 \xi_3=0,\\
x_5=0,\\
x_6 \xi_1+x_7\xi_2+x_8 \xi_3=0,\\
p_\rho x_1+(p_{S_2}-1)x_9=0,
\end{align*}
from which we obtain that there exist two linearly independent eigenvectors of the form $W=(0,x_2,x_3,x_4,0,0,0,0,0)^T$ with $(x_2,x_3,x_4)\cdot\xi=0$, two of the form $W=(0,0,0,0,0,x_6,x_7,x_8,0)^T$ with $(x_6,x_7,x_8)\cdot\xi=0$, and one of the form $W=(x_1,0,0,0,0,0,0,0,x_9)^T$ with $(x_1,x_9)\cdot (p_\rho,p_{S_2}-1)^T=0$.

Second, $g$ has four different zeros. Since $g$ is essentially the same as the corresponding one in one space dimension, we only (easily) transfer the considerations from \cite{HR3} for the reader's convenience.

Note that $g(\pm\infty)=+\infty$, and $g(0)=\left(\frac{\lambda (1-p_{S_2})}{\rho \tau_3}+p_\rho\right) \cdot \frac{\kappa}{\tau_1}>0$.

Let
\begin{align*}
\mu_\pm:=\pm \sqrt{\frac{\lambda(1-p_{S_2})}{\rho \tau_2}+p_\rho},
\end{align*}
then
\begin{align*}
\mu_-<0<\mu_+
\end{align*}
and
\begin{align*}
&g(\mu_\pm)=\mu_- \left( \frac{\theta p_\theta^2}{\rho^2 e_\theta} \mu_+ \mp\frac{\kappa \theta p_\theta}{\tau_1 \rho^2 e_\theta} p_q\cdot \xi \right)\equiv \mu_- Q,
\end{align*}
which implies by assumption \eqref{5.3} and \eqref{5.4} that
\begin{align*}
Q\ge \min_{(t,x)\in \mathbb R^+\times \mathbb R^3}\frac{ \theta  {p_\theta^2}  {\mu_+}}{2  {\rho^2}  {e_{\theta}}}>0
\end{align*}
if
\begin{align*}
|p_q\cdot \xi |\le \min_{(t,x)\in \mathbb R^+\times \mathbb R^3} \frac{\tau_1 p_\theta \mu_+}{\kappa} ,
\end{align*}
which is satisfied if $|q| <\delta$ for some $\delta>0$. Therefore,  we derive that
\begin{align*}
g( \mu_\pm)<0.
\end{align*}
Hence, $g$ has four different real zeros $z_1<z_2<0<z_3<z_4$.
Thus we have demonstrated that our system is\emph{ constantly }hyperbolic.

The constant hyperbolicity  implies the local well-posedness \cite{BS}. Constantly hyperbolic system are much less investigated than symmetric-hyperbolic or strictly hyperbolic ones, causing in general more difficulties. But we can refer to \cite[Thm. 2.3 and Thm. 10.2]{BS}, and thus have
\begin{theorem}\label{th5.1}
Let $s > \frac{n}{2} +1$ and $(\rho_0,u_0,\theta_0,q_0, S_{20}): \mathbb R^n\rightarrow \mathbb R^{2n+3}$ be given with
\begin{align*}
\rho_0-1, u_0, \theta_0-1, q_0, S_0 \in H^s,\\
\min_{x\in \mathbb R^n}(\rho_0(x), \theta_0(x))>0,\quad \max_{x\in \mathbb R^n} (|q_0(x)|, |S_{20}(x)|)\le \frac{\delta}{2}.
\end{align*}
Then, there exists a unique local solution $(\rho,u,\theta,q,S_2)$ to system \eqref{5.1} in some time interval $[0,T_0)$ with
\begin{align*}
(\rho-1,u,\theta-1,q,S_2) \in C^0([0,T_0),H^s) \cap C^1([0,T_0),H^{s-1})
\end{align*}
and
\begin{align*}
\min_{(x,t)\in \mathbb R^n \times [0,T_0)}(\rho(x,t), \theta(x,t))>0,\quad \max_{(x,t)\in  \mathbb R^n \times [0,T_0)} (|q(x,t)|, |S_{2}(x,t)|)\le \delta.
\end{align*}
\end{theorem}
 The following proposition states the finite propagation speed property which is guaranteed by the hyperbolicity of the system \eqref{5.1}, see \cite{Ra015, HR014}.
\begin{proposition}\label{pro5.1}
Assume the initial data $(\rho_0, u_0, \theta_0, q_0, S_{20})$ satisfy the assumption given in Theorem \ref{th5.1}  and $(\rho, u, \theta, q, S_2)$ be local solutions to \eqref{5.1} on $[0, T_0)$. We further assume that the initial data $(\rho_0-1, u_0, \theta_0-1, q_0, S_{20})$  are compactly supported in
a ball $B_0(M)$ with radius $M>0$. Then, there exists a constant $\sigma$ such that
\begin{align}\label{5.8}
(\rho(\cdot, t), u(\cdot, t), \theta(\cdot, t), q(\cdot, t), S_2(\cdot, t))=(1, 0, 1, 0, 0)=:(\bar \rho, \bar u, \bar \theta, \bar q, \bar S_2)
\end{align}
on $D(t):=\{x\in \mathbb R^n| |x|\ge M+\sigma t \},\, 0\le t < T_0$.
\end{proposition}
In the sequel, without loss of generality,  we shall assume that
\begin{align*}
\frac{\bar \theta}{2}<\theta<2 \bar \theta.
\end{align*}
We recall the averaged quantities defined in\eqref{5.9}, \eqref{5.10},
\begin{align*}
F(t)=\int_{\mathbb R^n}  x\cdot \rho(x,t) u(x,t)  \dif x,\\
G(t)=\int_{\mathbb R^n} ({\mathcal E}(x,t)-\bar {\mathcal E}) \dif x,
\end{align*}
with ${\mathcal E}(x,t)=\rho(e+\frac{1}{2} u^2)$ and $\bar {\mathcal E}=\bar \rho(\bar e+\frac{1}{2} \bar u ^2)=C_v$.  We mention that these quantities above exist as finite numbers since the solution $(\rho-1, u, \theta-1, q, S)$ is zero on the set $D(t)$ given in  Proposition \ref{pro5.1}.

Now, we are ready to prove Theorem \ref{th1.3}, with an ansatz and with calculations that we have been using starting in \cite{HR014}, essentially going back to ideas of Sideris \cite{Si84}, and that has also been used in the one-dimensional situation \cite{HuRaWa022}.
\begin{proof}
We will present the case $n=3$. The case $n=2$ is obtained with easy modifications. From the equations $\eqref{1.18}_{1-3}$, we can get the equation for the total energy $E$:
\begin{align} \label{5.11}
\partial_t E + \di (uE+up-uS_2+q)=0,
\end{align}
which implies that $G(t)$ is constant and
\begin{align}\label{New2}
G(t)=G(0)>0.
\end{align}
On the other hand, we have
\begin{align*}
F^\prime(t)&=\int_{\mathbb R^3} \partial_t(\rho u) x \dif x \\
&=\int_{\mathbb R^3} (-\di (\rho u \otimes u)-\nabla p +\nabla S_2 ) x \dif x \\
&=\int_{\mathbb R^3} \rho u^2 \dif x +3 \int_{\mathbb R^3} (p-\bar p) \dif x-3 \int_{\mathbb R^3} S_2 \dif x.
\end{align*}
By the constitutive equations \eqref{1.10} and \eqref{1.11}, we have
\begin{align*}
\int_{\mathbb R^3} (p-\bar p)\dif x = \int_{\mathbb R^3} (R \rho \theta- \frac{\tau_1}{2\kappa \theta} q^2-\frac{\tau_3}{2\lambda} S_2^2-R \bar \rho \bar \theta) \dif x
\end{align*}
and
\begin{align*}
R\rho \theta= \frac{R}{C_v} \rho e - \frac{\tau_1 R}{C_v \kappa \theta} q^2-\frac{\tau_3 R}{C_v \lambda} S_2^2.
\end{align*}
So, using \eqref{New2}, we derive that
\begin{align*}
\int_{\mathbb R^3} (p-\bar p) \dif x &\ge \int_{\mathbb R^3} \left( (\gamma-1) (\rho e- \bar \rho \bar e)-\frac{\tau_1 \gamma}{\kappa \theta} q^2-\frac{\tau_3 \gamma}{\lambda} S_2^2 \right) \dif x \\
&\ge \int_{\mathbb R^3} (\gamma-1)\left( (E-\frac{1}{2}\rho u^2)-\bar E\right) \dif x -\int_{\mathbb R^3} \left( \frac{\tau_1\gamma}{\kappa \theta} q^2+\frac{\tau_3 \gamma}{\lambda} S_2^2 \right)\dif x\\
&\ge -\frac{\gamma-1}{2} \int_{\mathbb R^3} \rho u^2 \dif x -\int_{\mathbb R^3} \left( \frac{\tau_1\gamma}{\kappa \theta} q^2+\frac{\tau_3 \gamma}{\lambda} S_2^2 \right)\dif x,
\end{align*}
where $\gamma:= \frac{R}{C_v}+1$. Combining the above estimates and using H\"{o}lder's inequality, we get
\begin{align} \label{5.12}
F^\prime (t)\ge \frac{5-3\gamma}{2} \int_{\mathbb R^3} \rho u^2 \dif x -3\int_{\mathbb R^3} \frac{\tau_1 \gamma}{\kappa \theta} q^2 \dif x
-3\int_{\mathbb R^3} \left(\frac{\tau_3 \gamma}{\lambda}+\frac{1}{2} \right) S_2^2 \dif x-2\pi(M+\sigma t)^3.
\end{align}
By the definition of $F(t)$, we conclude
\begin{align*}
F^2(t)&= \left( \int_{\mathbb R^3} x\cdot \rho(x,t) u(x,t) \dif x \right)^2 \\
&\le \int_{B_t} x^2 \rho \dif x \cdot \int_{B_t} \rho u^2 \dif x\\
&\le (M+\sigma t)^2 \int_{B_t} \rho \dif x \cdot \int_{B_t} \rho u^2 \dif x \\
&=(M+\tilde\sigma t)^2 \int_{B_0} \rho \dif x \cdot \int_{B_t} \rho u^2 \dif x \\
&\le \frac{4\pi}{3} \max \rho_0 (M+\tilde\sigma t)^5 \int_{\mathbb R^3} \rho u^2 \dif x,
\end{align*}
where $B_t:= \{ x\in \mathbb R^3| |x|\le M+\tilde \sigma t \}$ and $\tilde \sigma\ge \sigma$ can be chosen arbitrarily . For simplicity, we still denote
$\tilde \sigma$ by $\sigma$ in the following calculations. Therefore, we get
\begin{align} \label{5.13}
F^\prime(t) \ge \frac{3(5-3\gamma)}{8 \pi \max \rho_0 (M+\sigma t)^5 } F^2-3\int_{\mathbb R^3} \left( \frac{\tau_1 \gamma}{\kappa \theta} q^2+
\frac{2\tau_3\gamma+\lambda}{2\lambda} S_2^2 \right)\dif x-2\pi (M+\sigma t)^3.
\end{align}
Let $c_2:=\frac{\sigma}{M},\, c_3:= \frac{3 (5-3\gamma)}{8\pi \max \rho_0 M^5}$. Assume for the moment
\begin{align}\label{5.14}
F(t)\ge c_1>0
\end{align}
and
\begin{align} \label{5.15}
2\pi(M+\sigma t)^3 =2\pi M^3 (1+c_2t)^3 \le \frac{c_3}{2 (1+c_2t)^5} F^2,
\end{align}
where $c_1$ is to be determined later. Under the above a priori assumptions, we immediately get
\begin{align*}
F^\prime(t) \ge \frac{c_3}{2 (1+c_2 t)^5} F^2 -\frac{6\tau_1\gamma}{\kappa \bar \theta} \int_{\mathbb R^3} q^2 \dif x-\frac{6\tau_3\gamma+3\lambda}{2\lambda}
\int_{\mathbb R^3} S_2^2 \dif x.
\end{align*}
Using the assumption \eqref{5.14}, the above inequality implies that
\begin{align}\label{5.16}
\frac{F^\prime(t)}{F^2(t)}\ge \frac{c_3}{2 (1+c_2t)^5}-\frac{6\tau_1\gamma}{c_1^2 \kappa \bar \theta} \int_{\mathbb R^3} q^2 \dif x-\frac{6\tau_3\gamma+3\lambda}{c_1^2 2 \lambda}
\int_{\mathbb R^3} S_2^2 \dif x.
\end{align}
 Now, we recall the dissipative entropy equation \eqref{dissenteq} given in Lemma \ref{lem-dissenteq}, with $\mu=0$,
\begin{align}
 \partial_t\left[C_v \rho(\theta-\ln \theta -1)+R (\rho\ln \rho-\rho+1) +(1-\frac{1}{2\theta}) \frac{\tau_1}{\kappa \theta} q^2
 +\frac{1}{2}\rho u^2 +\frac{\tau_3}{2 \lambda} S_2^2\right] \nonumber\\
+ \di \left[\rho u C_v (\theta-\ln \theta -1)+u(1-\frac{1}{2\theta})\frac{\tau_1}{\kappa \theta} q^2+\frac{\tau_3}{2\lambda} u S_2^2
+R \rho u \ln \rho-R\rho u-\frac{q}{\theta}\right. \nonumber\\
\left.+\frac{1}{2}\rho u |u|^2+pu+q-S_2 u \right]+\frac{q^2}{\kappa \theta^2}+\frac{S_2^2}{\theta \lambda}=0. \label{New3}
\end{align}
Let
$$W_0=\int_{\mathbb R^3} \left(C_v \rho_0 (\theta_0-\ln \theta_0-1)+R(\rho_0 \ln \rho_0-\rho_0+1)+(1-\frac{1}{2\theta_0}) \frac{\tau_1}{\kappa \theta} q_0^2 +\frac{\tau_2}{2\lambda} S_{20}^2  \right) \dif x,$$
then \eqref{New3} implies
\begin{align*}
\int_0^t \int_{\mathbb R^3} \frac{q^2}{\kappa \theta^2} \dif x \dif t +\int_0^t \int_{\mathbb R^3} \frac{S_2^2}{\theta \lambda} \dif x \dif t \le W_0+\frac{\max \rho_0}{2} \|u_0\|_{L^2}^2.
\end{align*}
Therefore, we have
\begin{align}\label{5.17}
\frac{6\tau_1\gamma }{c_1^2 \kappa \bar \theta}\int_0^t \int_{\mathbb R^3} q^2 \dif x \dif t +\frac{6\tau_3 \gamma+3\lambda}{c_1^2 2 \lambda}\int_0^t \int_{\mathbb R^3} S_2^2 \dif x \dif t \le c_4+c_5 \|u_0\|_{L^2}^2,
\end{align}
where
\[
c_4=\frac{3}{c_1^2} \left[ \bar \theta (8\tau_1\gamma +2\tau_3 \gamma+\lambda)W_0\right], \quad c_5=\frac{3}{c_1^2} \left[ \bar \theta(8\tau_1\gamma +2\tau_3 \gamma+\lambda) \frac{\max\rho_0}{2}\right].
\]
Integrating \eqref{5.16} over $(0, t)$, using the above estimates, we have
\begin{align} \label{5.18}
\frac{1}{F_0}-\frac{1}{F(t)}\ge -\frac{c_3}{8c_2 (1+c_2 t)^4}+\frac{c_3}{8c_2}-c_4-c_5 \|u_0\|_{L^2}^2.
\end{align}
Now, we assume
\begin{align}
F_0 > \frac{16 c_2}{c_3}, \label{5.19}\\
c_4+c_5\|u_0\|_{L^2}^2 \le \frac{c_3}{16 c_2}. \label{5.20}
\end{align}
Then we have from \eqref{5.18} and \eqref{5.20}
\begin{align} \label{5.21}
\frac{1}{F_0} \ge \frac{1}{F_0}-\frac{1}{F(t)}\ge -\frac{c_3}{8c_2(1+c_2t)^4}+\frac{c_3}{16c_2}
\end{align}
which implies that the maximal time of existence $T$ can not arbitrarily large without contradicting \eqref{5.19}.

Now, we first show that the a priori assumption \eqref{5.14} holds. Define
$$
c_1:=\frac{4c_2}{c_3},
$$ then
\begin{align*}
\frac{1}{F(t)} \le \frac{1}{F_0} +\frac{c_3}{8c_2(1+c_2t)^4}-\frac{c_3}{16c_2}\le \frac{c_3}{8c_2(1+c_2t)^4}
\end{align*}
which implies
\begin{align} \label{5.22}
F(t)\ge \frac{8c_2}{c_3}(1+c_2t)^4\ge 2c_1.
\end{align}
This assures the a priori assumption \eqref{5.14} being compatible with \eqref{5.19} by noting that $F_0\ge 2c_1$.

To show the a priori estimates \eqref{5.15} hold, by bootstrap methods, we only need to show
\begin{align} \label{5.23}
2\pi M^3 (1+c_2t)^3\le \frac{c_3}{4(1+c_2t)^5}F(t)^2.
\end{align}
As a first step, we need \eqref{5.23} to hold for $t=0$, that is
\begin{align} \label{5.24}
F_0^2 \ge \frac{8\pi M^3}{c_3}=\frac{64 \pi \max \rho_0}{3(5-3\gamma)} M^8.
\end{align}
Using \eqref{5.22}, \eqref{5.23} is equivalent to
\begin{align} \label{5.25}
\sigma^2 \ge \frac{3(5-3\gamma)}{64 \max \rho_0}
\end{align}
which is satisfied naturally since $\sigma$ can be chosen arbitrarily large and then be fixed.

Thus, the proof will be finished if there exists $u_0$ such that \eqref{5.19}, \eqref{5.20} and \eqref{5.24} hold and the assumption \ref{1.16} is satisfied. Let (cp. \cite{HR014,HuRaWa022})
\begin{align} \label{5.26}
\tilde v(r)=
\begin{cases}
L \cos(\frac{\pi}{2}(r-1)),\quad\quad\quad\quad\quad  r\in [0,1],\\
L,\quad \quad\quad\quad \quad\quad\quad\quad\quad\quad\quad r\in(1, M-1],\\
\frac{L}{2} \cos ( \pi(r-M+1))+\frac{L}{2},\quad r\in (M-1, M],\\
0,\quad\quad\quad\quad\quad\quad\quad\quad\quad\quad\quad r\in(M, +\infty),
\end{cases}
\end{align}
where $L$ is a positive constant to be determined later. $\tilde v$ is not in $H^3(\mathbb{R}_+)$, but we can think of $\tilde v$ being smoothed around the singular points $r=1,M-1,M$ and put to zero around $r=0$, yielding a function $v$, with $\|v\|_{L^2}\leq 2 \|\tilde v\|_{L^2}$. We choose
$$
u_0(x):=v(|x|)\frac{x}{|x|}.
$$

 Assumption \ref{1.16} can easily be satisfied since it is equivalent to requiring
$$
\int_{\mathbb R^3}\left(\rho_0e_0 -\bar \rho \bar e + \frac{1}{2} u_0^2\right) dx > 0,
$$
which is satisfied by choosing $\rho_0\theta_0 >  \bar \rho \bar \theta=1$. Let $M\ge 5$, then
\begin{align*}
F_0&=\int_{\mathbb R^3} x \cdot \rho_0(x) u_0(x) \dif x =\int_{\mathbb R^3} \rho_0(x) v(|x|) |x|\dif x \\
&\ge \min \rho_0 \int_{B_0(M)} v(|x|) |x|\dif x \\
&\ge \min \rho_0 \int_0^M v(r) r\cdot 4\pi r^2 \dif r \\
&\ge \min \rho_0 \int _2^{M-2} L\cdot 4\pi r^3 \dif r\ge \frac{\pi \min \rho_0}{32} L M^4
\end{align*}
We choose $L$ sufficiently large such that
\begin{align*}
\frac{\pi \min \rho_0}{32} L \ge \max \left\{ \sqrt{\frac{64\pi \max \rho_0}{3(5-3\gamma)}}, \frac{128 \sigma \pi \max \rho_0}{3(5-3\gamma)}\right\}.
\end{align*}
So, \eqref{5.19} and \eqref{5.24} hold.
On the other hand, since $\|u_0\|_{L^2}^2 \le 4 L^2 \frac{4\pi}{3}M^3$, we choose $M$ such that
\begin{align*}
\bar \theta( 8\tau_1 \gamma+2\tau_2 \gamma+\mu)\left( W_0+\frac{2\pi\max \rho_0 L^2}{3}M^3\right)\le \frac{16 \pi \sigma \max \rho_0}{9(5-3\gamma)} M^4.
\end{align*}
 Therefore, \eqref{5.20} holds and the proof of Theorem \ref{th1.3} is finished.
\end{proof}


\section{Remark on the singular limit}

 We conclude the paper with a remark and additional result on the  \emph{singular limit $\tau\to 0$},
 which we can describe for the case $\tau_1>0, \tau_2=0, \tau_3>0$, $\mu>0$. We assume for simplicity $\tau_1=\tau_3=:\tau$.
For $\tau > 0$ let $(\rho^\tau,u^\tau,\theta^\tau, q^\tau, S_2^\tau)$ denote the  local solution to the system \eqref{1.8}
 defined on $[0, T_\tau)$,  where
\[
T_\tau=\sup \{T>0; (\rho^\tau-1,u^\tau,\theta^\tau-1, q^\tau, S_2^\tau)\in C([0,T],H^3),(\rho^\tau, u^\tau, \theta^\tau, q^\tau, S_2^\tau)\in G_1\}
\]
with initial data $(\rho_0^\tau, u_0^\tau, \theta_0^\tau, q_0^\tau, S_{20}^\tau)\in G_0$.
Then we have
\begin{theorem}\label{th1.2}
Let $(\rho,u,\theta)$ be the smooth solution to the classical compressible Navier-Stokes equations with $(\rho(x,0),u(x,0),\theta(x,0))=(\rho_0,u_0, \theta_0)$ satisfying
$\inf_{(x,t)\in \mathbb R^3 \times[0,T_*]} (\rho(x,t), \theta(x,t)) >0$ and
\begin{align*}
(\rho-1)\in C([0,T_*],H^6)\cap C^1([0,T_*],H^5),\\
(u, \theta-1)\in C([0,T_*],H^6) \cap C^1([0,T_*],H^4),
\end{align*}
with finite $T_*>0$. Moreover,  assume that the initial data are well-prepared, i.e.,
\begin{align*}
\|(\rho_0^\tau -\rho_0, u_0^\tau- u_0, \theta_0^\tau-\theta_0, \sqrt{\tau}(q_0^\tau+\kappa \nabla \theta_0), \sqrt{\tau}(S_{20}^\tau- \lambda \di u_0))\|_{H^3} \le  \tau.
\end{align*}
Then, there exist constants $\tau_0$ and $C>0$ such that for $\tau\le \tau_0$,
\begin{align}\label{5n.1}
\|(\rho^\tau,u^\tau,\theta^\tau)(\cdot,t)-(\rho,u, \theta)(\cdot,t)\|_{H^3}\le C \tau,
\end{align}
and
\begin{align}\label{5n.2}
\|(q^\tau+\kappa \nabla \theta, S_2^\tau-\lambda \nabla u) \|_{H^3} \le C \tau^\frac{1}{2},
\end{align}
for all $t\in (0,min(T_*,T_\tau))$, and the constant $C$ is independent of $\tau$.
\end{theorem}
Thus, we have that the solutions for $\tau>0$ converge to the solution for $\tau=0$ with a certain order in $\tau$ on any finite interval of common existence.

The (long) proof of Theorem \ref{th1.2}, which we omit here, can be done in the spirit of corresponding considerations in \cite{HR2}, overcoming a higher complexity given here by energy estimates similar to those used in the proof of Theorem \ref{th1.1}.

Recently, Peng and Zhao \cite{PeZh022} studied the 1-d version and obtained in partical a global existence result which is uniform with respect to $\tau$ as well as a global convergence result in a weak topology.


\begin{thebibliography}{aaaaa}

\bibitem{BS} S. Benzoni-Gavage and D. Serre, {\it Multi-dimensional hyperbolic partial differential equations: first-order system and applications,} Clarendon Press, Oxford (2007).
\bibitem{CA}  C. Cattaneo, Sulla coduzione del calore,  {\it Atti Sem. Mat. Fis. Univ. Modena} {\bf 3} (1948), 83-101.

\bibitem{CS} D. Chakraborty and J.E. Sader, Constitutive models for linear compressible viscoelastic flows of simple liquids at nanometer length scales,
{\it Physics of Fluids} {\bf 27} (2015), 052002-1--052002-13.

\bibitem{CG} P.J. Chen and M.E. Gurtin, On second sound in materials with memory, {\it Z. Ang. Math. Phys.} {\bf 21} (1970), 232-241.

\bibitem{CJ}{Y. Cho and B.J. Jin, Blow-up of viscous heat-conducting compressible flows, {\it J. Math. Anal. Appl.} {\bf 320} (2) (2006), 819-826. }

 \bibitem{CJo} C.I. Christov and P.M. Jordan, Heat condction paradox involving second-sound propagation in moving media, {\it Phys. Rev. Letters} {\bf 94} (2005), 154301-1---154301-4.

\bibitem{CFO} B.D. Coleman, M. Fabrizio and D.R. Owen, On the thermodynamics of second sound in dielectric crystals, {\it Arch. Rational Mech. Anal.} {\bf 80} (1986), 135-158.

\bibitem{CHO} B.D. Coleman, W.J. Hrusa and D.R. Owen, Stability of Equilibrium for a Nonlinear Hyperbolic System Describing Heat Propagation by Second Sound in Solids, {\it Arch. Rational Mech. Anal.} {\bf 94} (1986), 267-289.

\bibitem{FE}  E. Feireisl, A. Novotny and H. Petzeltov\'{a}, On the existence of globally defined weak solutions to
the Navier-Stokes equations, {\it J. Math. Fluid Mech.} {\bf 3} (2001), 358-392.

\bibitem{FeMu012} H.D. Fern\'andez Sare and J.E. Mu\~noz Rivera, Optimal
rates of decay in 2-d thermoelasticity with second sound,
{\it J. Math. Phys.} {\bf 53} (2012), 073509.

 
\bibitem{DA1} D. Hoff, Global existence for 1D, compressible, isentropic Navier-Stokes equations
with large initial data, {\it Trans. Amer. Math. Soc.} {\bf 303} (1) (1987), 169-181.

\bibitem{DA2} D. Hoff, Global solutions of the Navier-Stokes equations for multidimensional compressible
flow with discontinuous initial data, {\it J. Differential Equations} {\bf 120} (1) (1995), 215-254.


\bibitem{HW} Y. Hu and N. Wang, Global existence versus blow-up results for one dimensional compressible Navier-Stokes equations with Maxwell's law, {\it Math. Nachr.} {\bf 292} (2019), 826-840.

\bibitem{HR014} Y. Hu and R. Racke, Formation of singularities in one-dimensional
   thermoelasticity with second sound. {\it Quart. Appl. Math.} {\bf 72} (2014), 311-321.

\bibitem{HR1} Y. Hu and R. Racke, Compressible Navier-Stokes equations with hyperbolic heat
conduction, {\it J. Hyper. Diff. Equ.} {\bf 13} (2) (2016), 233-247.

\bibitem{HR2} Y. Hu and R. Racke, Compressible Navier-Stokes equations with revised Maxwell's law, {\it J. Math.Fluid Mech.} {\bf 19} (2017), 77-90.

\bibitem{HR3} Y. Hu and R. Racke, Hyperbolic compressible Navier-Stokes equations, {\it J. Differential Equations} {\bf 269} (2020), 3196-3220.

\bibitem{HuRaWa022}  Y. Hu, R. Racke and N. Wang, Formation of singularities for one-dimensional relaxed compressible Navier-Stokes equations, {\it J. Differential Equations} (accepted) (2022).

 \bibitem{JR00}  S. Jiang and R. Racke, {\it Evolution equations in thermoelasticity.}
   $\pi$ Monographs Surveys Pure Appl. Math.
  {\bf 112}. Chapman \& Hall/CRC, Boca Raton (2000).

  \bibitem{JZ01} S. Jiang and P. Zhang, Global spherically symmetry solutions of the compressible isentropic
Navier-Stokes equations, {\it Comm. Math. Phys.} {\bf 215} (2001), 559-581.
\bibitem{JZ03} S. Jiang and P. Zhang, Axisymmetric solutions of the 3-D Navier-Stokes equations for compressible isentropic fluids, {\it J. Math. Pures. Appl.} {\bf 82} (2003), 949-973.
 
\bibitem{KW} S. Kawashima, {\it Systems of a hyperbolic-parabolic composite type, with applications to the equations of magnetohydrodynamics,} Thesis, Kyoto University (1983).

 \bibitem{LI} P.L. Lions, {\it Mathematical Topics in Fluid Mechanics}, Vol.II, Compressible Models. Clarendon Press, Oxford (1998).

\bibitem{Ma84} A. Majda, {\it Compressible fluid flow and systems of
  conservation laws in several space variables.} Appl.\ Math.\ Sci.\ {\bf 53}, Springer, New York (1984).

  \bibitem{MN} A. Matsumura and T. Nishida, The initial value problem for the
equations of motion of viscous and heat-conductive gases, {\it J. Math.
Kyoto Univ.} {\bf 20} (1) (1980), 67-104.

 
\bibitem{MA} J.C. Maxwell, On the dynamics theory of gases, {\it Phil. Trans. R. Soc. Lond.} {\bf 157} (1867), 49-88.

\bibitem{NA} J. Nash, Le probl\`{e}me de Cauchy pour les \'{e}quations
diff\'{e}rentielles d'un fluide g\'{e}n\'{e}ral, {\it Bull. Soc. Math.
France} {\bf 90} (1962), 487-497.

\bibitem{PC} M. Pelton, D. Chakraborty,E. Malachosky, P. Guyot-Sionnest and J. E. Sader, Viscoelastic flows in simple liquids generated by vibrating nanostructures,
{\it Phys. Rev. Lett.}, {\bf 111} (2013), 244502.

\bibitem{PeZh022} Y.-J. Peng and L. Zhao, Global convergence to compressible full Navier–Stokes equations by approximation with Oldroyd-type constitutive laws, {\it J. Math. Fluid Mech.} {\bf (2022)},  24:29.
 
\bibitem{Ra015} R. Racke, {\it Lectures on Nonlinear Evolution Equations. Initial Value Problems,} 2nd edition, Birkh\"auser, Basel (2015).

 
\bibitem{SE} J. Serrin, On the uniqueness of compressible fluid motion, {\it Arch. Rational Mech. Anal.} {\bf 3} (1959), 271-288.

\bibitem{ShKu019} B. Sharma and R. Kumar, Estimation of bulk viscosity of dilute gases using a nonequlibirium molecular dynamics approach. {\it Pys. Rev. E {\bf 100}} (2019), 013309.

\bibitem{SK85} Y. Shizuta and S. Kawashima, Systems of equations of hyperbolic-parabolic type with applications to the discrete Boltzmann equation, {\it Hokkaido Math. J.} {\bf 14} (1985), 249-275.

\bibitem{Si84} T.C. Sideris, Formation of singularities in solutions to nonlinear hyperbolic equations, {\it Arch. Rational Mech. Anal.} {\bf 86} (1984), 36909000938

 \bibitem{TA}{M.A. Tarabek, On the existence of smooth solutions in one-dimensional
nonlinear thermoelasticity with second sound, {\it Quart. Appl. Math.} {\bf 50}
 (1992), 727--742.}

 \bibitem{UmKaSh84} T. Umeda, S. Kawashima and Y. Shizuta, On the decay of solutions to the linearized equations of electro-magneto-fluid dynamics, Japan J. Appl. Math. {\bf 1} (1984), 435-457.


\bibitem{X} Z.P. Xin, Blowup of smooth solutions to the compressible Navier-Stokes equation
with compact density, {\it Comm. Pure. Appl. Math.} {\bf 51} (1998), 229-240.

 \bibitem{YW14} W.A. Yong, Newtonian limit of Maxwell fluid flows, {\it Arch. Rational Mech. Anal.} {\bf 214} (2014), 913-922.



\end{thebibliography}
\end{document}